\documentclass[12pt,reqno]{amsart}
\usepackage{amssymb,amsfonts,amsthm,amsmath,mathrsfs}
\usepackage{cite}
\usepackage[shortlabels]{enumitem}
\usepackage[left=1 in,top=1 in,right=1 in, bottom=1 in]{geometry}
\usepackage{graphicx}
\usepackage{color}
\usepackage[pagebackref=false]{hyperref}
\usepackage{tcolorbox}
\newcommand{\sign}{{\rm sign}}

\def\d{{\rm d}}
\def\eqdef{\stackrel{\rm def}{=}}
\definecolor{darkred}{rgb}{.70,.12,.20}

\definecolor{darkgreen}{rgb}{.20,.52,.14}

\definecolor{byz}{rgb}{.44,.16,.39}

\numberwithin{equation}{section}

\newtheorem{theorem}{Theorem}[section]
\newtheorem{lemma}[theorem]{Lemma}

\newtheorem{corollary}[theorem]{Corollary}
\newtheorem{assumption}[theorem]{Assumption}

\newtheorem{proposition}[theorem]{Proposition}

\theoremstyle{remark}
\newtheorem{remark}[theorem]{Remark}
\newtheorem{example}[theorem]{\bf{Example}}

\usepackage{todonotes}

\newcommand{\varep}{\varepsilon}

\newcommand{\beq}{\begin{equation}}
\newcommand{\eeq}{\end{equation}}
\newcommand{\beqs}{\begin{equation*}}
\newcommand{\eeqs}{\end{equation*}}
\newcommand{\ba}{\begin{array}}
\newcommand{\ea}{\end{array}}
\newcommand{\beas}{\begin{eqnarray*}}
\newcommand{\eeas}{\end{eqnarray*}}
\newcommand{\bea}{\begin{eqnarray}}
\newcommand{\eea}{\end{eqnarray}}
\newcommand{\bal}{\begin{align}}
\newcommand{\eal}{\end{align}}

\newcommand{\bals}{\begin{align*}}
\newcommand{\eals}{\end{align*}}

\newcommand{\tnum}{\rm(\roman*)}
\newcommand{\rnum}{\rm(\alph*)}

\newcommand{\R}{\ensuremath{\mathbb R}}

\newcommand{\N}{\ensuremath{\mathbb N}}

\newcommand{\bds}{\begin{displaystyle}}
\newcommand{\eds}{\end{displaystyle}}

\newcommand{\remove}[1]{} 
\renewcommand{\remove}[1]{#1} 

\definecolor{darkred}{rgb}{.70,.12,.20}

\definecolor{darkgreen}{rgb}{.20,.52,.14}

\title{
Asymptotic estimates for solutions of inhomogeneous non-divergence diffusion equations with drifts
}

\author{Luan Hoang$^{1,*}$}
\address{$^{1}$Department of Mathematics and Statistics,
Texas Tech University,
1108 Memorial Circle, Lubbock, TX 79409--1042, U. S. A.}
\email{luan.hoang@ttu.edu}

\author{Akif Ibragimov$^{1,2}$}
\address{$^{2}$Oil and Gas Research Institute, Russian Academy of Science,  3 Gubkin Street, Moscow, 119333, Russia}
\email{ilya1sergey@gmail.com}

\thanks{$^*$Corresponding author.}

\date{\today}


 
\begin{document}

\begin{abstract} 
We study  the long-time dynamics of the nonlinear processes modeled by diffusion-transport partial differential equations in non-divergence form with drifts. The solutions are subject to some inhomogeneous Dirichlet boundary condition.  Starting with the reduced linear problem, we obtain the asymptotic estimates for the solutions, as time $t\to\infty$, depending on the asymptotic  behavior of the forcing term and boundary data. These are established in both cases when the drifts are uniformly bounded, and unbounded as $t\to\infty$. For the nonlinear problem, we prove the convergence of the solutions under suitable conditions that balance the  growth of the nonlinear term with the decay of the data. To take advantage of the diffusion in the non-divergence form, we prove an inhomogeneous version of the Landis-typed  Growth Lemma and apply  it to successive time-intervals. At each time step, the center for the barrier function is selected carefully to optimize the contracting factor. Our rigorous results show the robustness of the model.
\end{abstract}

\maketitle 
\tableofcontents 

\pagestyle{myheadings}\markboth{\sc L.~Hoang and A.~Ibragimov}
{\sc Inhomogeneous non-divergence diffusion equations with drifts}

\section{Introduction }\label{intro}
We consider a fluid flow in porous media in the $n$-dimensional Euclidean space. (For physical problems, $n\le 3$). 
Let $x\in \mathbb R^n$ be the spatial variable and  $t\in\R$ the time variable.
Let $\rho(x,t)$, $p(x,t)$ and $v(x,t)$ denote  the density, pressure and velocity  of the fluid at location $x$ and time $t$. 
In \cite{HI3}, we derived an equation based on Eistein's probabilistic method for the Brownian motion \cite{Einstein1905}, namely,
\beq\label{drveq}
\frac{\partial \rho}{\partial t}=\langle A(x,t),D^2\rho\rangle  + (M_0(x,t) v(x,t))\cdot \nabla \rho,
\eeq
where $A(x,t)$ and $M_0(x,t)$ are  $n\times n$ matrices such that
\beqs 
\xi^{\rm T} A(x,t)\xi > 0 \text{ and }\xi^{\rm T} M_0(x,t)\xi\ge 0\text{ for all }\xi\in\R^n.
\eeqs 
In equation \eqref{drveq}, the term $\langle A(x,t),D^2\rho\rangle$ denotes the inner product between the matrix $A(x,t)$ and the Hessian matrix $D^2\rho$, see the notation in section \ref{maxmax} below. It represents the diffusion in the non-divergence form, while the term $(M_0(x,t) v(x,t))\cdot \nabla \rho$ represents the transport.

Assume the anisotropic Darcy's law \cite{DarcyBook,BearBook}, 
\beq \label{Darcy}
v=-\bar K(x,t)(\nabla p - \rho \vec g), 
\eeq 
where $\bar K(x,t)$ is an $n\times n$ matrix, and $\vec g$  is the gravitational acceleration in the case $n=1,2,3$ or any constant vector otherwise.
Combining \eqref{drveq} with \eqref{Darcy} yields
\beq\label{eqgrav}
\frac{\partial \rho}{\partial t}=\langle A(x,t),D^2\rho\rangle -(K_0 (x,t)\nabla p)\cdot \nabla \rho + \rho B_0(x,t)\cdot \nabla \rho,
\eeq
where
\beqs 
 K_0(x,t)=M_0(x,t)\bar K(x,t),\quad B_0(x,t)=M_0(x,t) \bar K(x,t)\vec g.
\eeqs


When additional relations between $\rho$ and $p$ are specified, we can derive from \eqref{eqgrav} a nonlinear partial differential equation (PDE) for $\rho$.
The following are some typical cases in fluid mechanics.

\medskip\noindent
\textit{Case of isentropic gas flows.} In this case,  $p=c\rho^\gamma$ for some  constant $c>0$, where  $\gamma\ge 1$ is the specific heat ratio. Then equation \eqref{eqgrav} has the specific form
\beq\label{iseneq}
\frac{\partial \rho}{\partial t}=\langle A(x,t),D^2\rho\rangle - c\gamma \rho^{\gamma-1} (K_0(x,t)\nabla \rho)\cdot \nabla \rho +  \rho B_0(x,t)\cdot \nabla \rho.
\eeq
In particularly,  one has $\gamma=1$ for ideal gases and equation \eqref{iseneq} reads as
\beq\label{idealeq}
\frac{\partial \rho}{\partial t}=\langle A(x,t),D^2\rho\rangle - c(K_0(x,t)\nabla \rho)\cdot \nabla \rho +  \rho B_0(x,t)\cdot \nabla \rho.
\eeq

\medskip\noindent
\textit{Case of slightly compressible fluids.} The equation of state is 
$$\frac1\rho \frac{d\rho}{dp}=\kappa, \text{ where $\kappa$ is the small, positive, constant compressibility.}$$
Then \eqref{eqgrav} can be rewritten as 
\beq\label{slighteq}
\frac{\partial \rho}{\partial t}=\langle A(x,t),D^2\rho\rangle  - \frac{1}{\kappa\rho} (K_0(x,t)\nabla \rho)\cdot \nabla \rho + \rho B_0(x,t) \cdot \nabla \rho.
\eeq

Equations \eqref{iseneq}, \eqref{idealeq} and \eqref{slighteq} fall into a larger class of equations which is formulated by \eqref{maineq} below. This nonlinear parabolic equation in non-divergence form will be the subject of our study. 
To the best of our knowledge, there are not many papers dealing with the long-time dynamics of non-divergence parabolic equations, especially the highly nonlinear ones such as \eqref{maineq}. 
See \cite{KrylovBookHolder,KrylovBookSobolev} for the general existence and regularity theory.
See also \cite{Safonov2010} for results for non-divergence linear elliptic problems with unbounded drift.
In our previous work \cite{HI3}, we essentially study \eqref{maineq} with zero forcing and, after shifting the solution by a constant, the homogeneous boundary condition. The purpose of the current paper is to obtain results for inhomogeneous problems. We emphasize that the equations we are dealing with contain drifts which can be nonlinear and unbounded as time tends to infinity. As we will see, it requires more technical and sophisticated tools.

There are three ingredients that we develop in the current paper to study of this highly nonlinear problem. The first is the general Bernstein--Cole--Hopf transformations that convert the equation to more manageable inequalities. The second is the Growth Lemma for inhomogeneous problems expressed in terms of the time-integral of the forcing function. The third is the iteration of this lemma on successive time-intervals with appropriate adjustments at each step. With those, we are able to establish the convergence/decay of the solutions as time tends to infinity. 

The paper is organized as follow.
In section \ref{prelim}, we recall the notation, set the background and  introduce the main initial boundary value problem \eqref{nonIBVP}.
 In section \ref{maxmax}, we use the standard Maximum Principle to derive preliminary bounds for the solutions and sub-/super-solutions of inhomogeneous problems for both linear and nonlinear equations. A new feature compared to the usual form in literature is  that the bounds are expressed in terms of the integral in time for the forcing function. This hints at  estimates that we will obtain and explore later which are similar to the ones for the two-dimensional Navier--Stokes equations by Foias and Prodi \cite{FP67}.
 In section \ref{linsec}, we study the linear problem. We recall in  Lemma \ref{lemgrowth} the Growth Lemma in a finite interval for the homogeneous problem.
 This is extended to the inhomogeneous problem in Lemma \ref{lemG2}.
 A general but complicated estimate for all time is established in Proposition \ref{genlem}. The main idea is combining the application of the Growth Lemma in successive intervals of time with the moving center in the construction of the barrier function at each step. This allows us to gain better contracting factors of the Growth Lemma (not too close to $1$.)
 The estimate in Proposition \ref{genlem} is then used to establish a particular decay estimate for the solution in the case of bounded drifts, see Theorem \ref{NHthm}.
 The case of unbounded drifts as $t\to\infty$ is treated in Theorems \ref{NHthm3} and  \ref{NHthm2}. Basically, even if the drift term goes to infinity as $t\to\infty$, we still have the boundedness or even the decay for the solution if this growth is not too fast and the forcing function and boundary data decay at an appropriate rate to balance out that growth.
Although the assumptions are technical, some explicit and simple examples are provided in Propositions \ref{ex1} and \ref{ex2} to demonstrate the applicability of the results. 
In section \ref{nonlinsec}, we study the nonlinear problem. 
We recall from our previous work \cite{HI3}  the transformation of Bernstein--Cole-Hopf type \cite{Bernstein1912,Cole1951,Hopf1950} in order to remove the quadratic term of the gradient. 
However, unlike \cite{HI3} the resulting partial differential inequalities \eqref{Lw1} and \eqref{Lw2} are still highly nonlinear. Nonetheless, they will be treated by the techniques developed in the preceding section for linear equations with possible unbounded drifts at time infinity.
The case for $L^1_t$-forcing function is treated in Theorems \ref{NHL1} and \ref{NHL2}.
Note that even when the solution is proved to be bounded, the contributing drift term still may not be so, see \eqref{B3}.
For the case when the forcing function is not in class $L_t^1$, we are able to treat the slightly compressible fluids in Theorems \ref{NHL3} and \ref{NHL4}. 

It is worth mentioning that the obtained results show the robustness of the dynamical system generated by the studied PDE. This in part justifies our models \eqref{drveq}--\eqref{slighteq} of the filtration for a class of fluid flows in porous media.
 
\section{Preliminaries}\label{prelim}

Throughout the paper, the spatial dimension $n\ge 1$ is fixed. 
For a vector $x\in\R^n$, its Euclidean norm is denoted by $|x|$.
 Let $\mathcal M^{n\times n}$ denote the set of $n\times n$ matrices of real numbers, and $\mathcal M^{n\times n}_{{\rm sym}}$ denote the set of symmetric matrices in $\mathcal M^{n\times n}$.
For two matrices $A,B\in \mathcal M^{n\times n}$, their inner product $\langle A,B\rangle$ is the trace ${\rm Tr}(A^{\rm T}B)$. For a scalar function $f(x)$ with $x=(x_1,\ldots,x_n)\in\R^n$, we denote by $D^2f$ the Hessian matrix of the second partial derivatives $(\partial^2f/\partial x_i\partial x_j)_{i,j=1,\ldots,n}$.

For a real-valued function $f$, we recall that the positive and negative parts of $f$ are  $f^+=\max\{0,f\}$ and $f^-=\max\{0,-f\}=(-f)^+$. Then one has    
\begin{align*}
f=f^+-f^-,\  |f|=f^+ + f^-=\max\{f^+,f^-\},\\ 
0\le f^+\le |f|,\  0\le f^-\le |f|,
\text{ and } 
-f^-\le f\le f^+.
\end{align*} 
If a function $f:S\to \R$ is bounded from above, then
$\displaystyle \max\Big\{0,\sup_S f\Big\}=\sup_S (f^+)$.

For the remaining of the paper, $U$ is a non-empty, open, bounded subset of  $\R^n$ with boundary $\Gamma$ and closure $\bar U$. Denote 
$d={\rm diam} (U)=\max\{|x-x'|:x,x'\in \bar U\}.$
For any point $y\not\in \bar U$, set
\beqs
r_*(y)=\min\{|x-y
|:x\in\bar U\}\text{ and } R_*(y)=\max\{|x-y|:x\in\bar U\}.
\eeqs
Then one has, for any $y\not\in \bar U$, 
\beq \label{rRd}
R_*(y)-d\le r_*(y)<R_*(y)\le r_*(y)+d.
\eeq 

We denote by $C_{x,t}^{2,1}(U\times I)$, for any interval $I$ of $\R$, the class of continuous functions $u(x,t)$ on $U\times I$ with continuous partial derivatives $\partial u/\partial t$, $\partial u/\partial x_i$, $\partial^2 u/\partial x_i\partial x_j$ for $i,j=1.\ldots,n$.

For a number $T\in(0,\infty)$, 
denote $U_T=U\times (0,T]$ and the parabolic boundary $\Gamma_T=\overline{U_T}\setminus U_T$, where $\overline{U_T}=\bar U\times[0,T]$ is the closure of $U_T$ (in $\R^{n+1}$). 
For the sake of convenience, we also denote $U_\infty=U\times(0,\infty)$ and its boundary $\Gamma_\infty=\partial (U_\infty)$.

Hereafter, $J$ is an interval  in $\R$ with non-empty interior, and $P$ is a function in $C^1(J,\R)$ with the derivative 
\beq \label{Pd}
P'\in C(J,[0,\infty)).
\eeq 

For fluid flows in porous media, $P$ is related to the equation of state for $\rho$ and $p$. However, we will consider a general function $P$.
It is clear from \eqref{Pd} that $P$ is an increasing function on $J$.

Let an extended number $T\in(0,\infty]$ be fixed. 
Let $A:U_T\to \mathcal M^{n\times n}_{{\rm sym}}$, $K:U_T\to \mathcal M^{n\times n}$ and $B:U_T\times J\to \R^n$  be given functions.

We  study the nonlinear equations \eqref{iseneq}, \eqref{idealeq} and \eqref{slighteq} in the general form
\beq \label{maineq}
\frac{\partial u}{\partial t}-\langle A(x,t),D^2u\rangle +   B(x,t,u)\cdot \nabla u+ P'(u) (K(x,t)\nabla u)\cdot \nabla u=f(x,t),
\eeq 
where $f:U_T\to \R$ is a given function.

For example, for slightly compressible fluids and equation \eqref{slighteq}, we can use
\beq \label{slighchoice}
J=(0,\infty),\quad P(s)=\ln s,\quad K(x,t)=\kappa^{-1}K_0(x,t),\quad B(x,t,s)=-sB_0(x,t).
\eeq 
For isentropic, non-ideal gas flows and equation \eqref{iseneq} with $\gamma>1$, we can use
\beq \label{isenchoice}
J=[0,\infty),\quad P(s)=s^\gamma,\quad K(x,t)=cK_0(x,t),\quad B(x,t,s)=-sB_0(x,t).
\eeq 
For ideal gas flows and equation \eqref{idealeq}, we can still use \eqref{isenchoice} with $\gamma=1$ for physical $u=\rho$ (the density), but can also alter $J=\R$ in \eqref{isenchoice} for mathematical $u$. Thus,
\beq\label{idealchoice} 
J=[0,\infty) \text{ or } J=\R,\quad  P(s)=s,\quad K(x,t)=cK_0(x,t),\quad B(x,t,s)=-sB_0(x,t).
\eeq

\begin{remark}
Although we focus on equation \eqref{maineq} with the constraint \eqref{Pd} in this paper, our obtained results below will apply to a more general class of equations, namely,
\beq \label{Lmu}
\frac{\partial u}{\partial t}-\langle A(x,t),D^2u\rangle +   B(x,t,u)\cdot \nabla u+ \mu(u) (K(x,t)\nabla u)\cdot \nabla u=f(x,t),
\eeq 
where $\mu$ is a continuous function from $J$ to $[0,\infty)$ or to $(-\infty,0]$.
Indeed, in the former case, we just take $P$ such that $P'=\mu$, while in the latter case, we take $P'=-\mu$ and  the new $K:=-K$. Then equation \eqref{Lmu}  is of the form \eqref{maineq} with the function $P$ satisfying \eqref{Pd}.
\end{remark}

We define the nonlinear operator $L$ associated with the left-hand side of \eqref{maineq} explicitly as
\beq \label{Ldef}
Lu=\frac{\partial u}{\partial t}-\langle A(x,t),D^2u\rangle +   B(x,t,u)\cdot \nabla u+ P'(u) (K(x,t)\nabla u)\cdot \nabla u
\eeq 
for any function $u\in C_{x,t}^{2,1}(U_T)$ with $u(U_T)\subset J$. 
The main task of this paper is to study the following initial boundary value problem
\beq\label{nonIBVP}
\begin{cases}
Lu(x,t)=f(x,t)&\text{ in }U\times(0,\infty),\\
    u(x,t)=g(x,t)&\text{ on } \Gamma\times(0,\infty),\\
    u(x,0)=u_0(x)&\text{ in }U,
\end{cases}
\eeq
where the forcing term $f:U_\infty\to \R$, the boundary data $g:\Gamma\times(0,\infty)\to \R$ and the initial condition $u_0:U\to \R$ are given functions.

A  function $u(x,t)$ is said to be a solution of Problem \eqref{nonIBVP} if  it belongs to $C_{x,t}^{2,1}(U_\infty)$ with  $u(U_\infty)\subset J$
and satisfies all equations in \eqref{nonIBVP}.

In order to study \eqref{nonIBVP}, we need some fine properties of solutions or sub-/super-solutions of a related linear inhomogeneous, non-divergence equations with general drifts. It involves the general linear operator of the following form. Given a function $b:U_T\to \R^n$, we define the linear operator $\widetilde L $ by
\beq \label{Ltil}
\widetilde L w = w_t-\langle A(x,t),D^2w\rangle  +b(x,t)\cdot \nabla w \text{ for }w\in C_{x,t}^{2,1}(U_T).
\eeq

\section{Direct estimates by the Maximum Principle}\label{maxmax}
Below, we establish a maximum principle for solutions of the inhomogeneous problem
$\mathscr Lu(x,t)=f(x,t)$ with either nonlinear operator $\mathscr L=L$ in \eqref{Ldef} or the linear one $\mathscr L=\widetilde L$ in \eqref{Ltil}.

\begin{theorem}\label{maxprin2}
Given a number $T\in(0,\infty)$, assume that there exists a constant $c_0>0$ such that
\beq\label{Aelip}
\xi^{\rm T} A(x,t)\xi \ge c_0|\xi|^2\text{ for all $(x,t)\in U_T$ and all $\xi\in \R^n$.}
\eeq 
Let $u\in C(\overline{U_T})\cap C_{x,t}^{2,1}(U_T)$ and define, for $t\in(0,T]$,
\beqs
G_1(t)=\max_{\Gamma_t} u, \quad G_2(t)=\min_{\Gamma_t} u,\quad G(t)=\max_{\Gamma_t} |u|.
\eeqs
Consider  the operator $\mathscr L$ to be either
\begin{enumerate}[label=\rnum]
   \item\label{ma}  $\mathscr L=\widetilde L$, or
   \item\label{mb} $\mathscr L=L$ with $u(U_T)\subset J$.
\end{enumerate} 
Then the following statements hold true.
\begin{enumerate}[label=\tnum]
\item\label{MP1}  Suppose there is a functions $F_1\in C([0,T],[0,\infty))$ such that 
    $\mathscr Lu(x,t)\le F_1(t)$ for all $(x,t)\in U_T$. 
Then
    \beq\label{max2}
        \max_{\overline{U_t}} u \le G_1(t)+\int_0^t F_1(\tau)\d \tau \text{ for all }t\in(0,T].
    \eeq

\item\label{MP2}  Suppose there is a  function $F_2\in C([0,T],[0,\infty))$  such that 
    $\mathscr Lu(x,t)\ge -F_2(t)$ for all $(x,t)\in U_T$. 
Then
\beq\label{min2}
        \min_{\overline{U_t}} u \ge G_2(t)-\int_0^t F_2(\tau)\d \tau \text{ for all }t\in(0,T].
    \eeq

\item\label{MP3} Suppose there is a  function $F\in C([0,T],[0,\infty))$  such that
       $|\mathscr Lu(x,t)|\le F(t)$ for all $(x,t)\in U_T$. 
Then
    \beq\label{maxabs}
        \max_{\overline{U_t}} |u| \le G(t)+\int_0^t F(\tau)\d \tau \text{ for all }t\in(0,T].
        \eeq
\end{enumerate}
\end{theorem}
\begin{proof}
We consider Case \ref{ma} first.

\medskip
\ref{MP1}  Define, for $t\in[0,T]$, 
$\mathcal{F}_1(t)=\int_0^t F_1(t)\d \tau.$
 Set $v(x,t)=u(x,t)- \mathcal{F}_1(t)$. Then
        $$\widetilde L v(x,t)=\widetilde  Lu(x,t) -\mathcal{F}_1'(t)=Lu(x,t)- F_1(t)\le 0.$$
    On $\Gamma_t$, one has $v(x,t)\le u(x,t)\le G_1(t)$. By the Maximum Principle for $\widetilde L$ and $v$, we have  $v(x,t)\le G_1(t)$ on $U_T$. Thus, $u(x,t)\le G_1(t)+\mathcal{F}_1(t)$ on $U_T$, which proves \eqref{max2}.

\medskip
 \ref{MP2} Define, for $t\in[0,T]$, $\mathcal{F}_2(t)=\int_0^t F_2(t)\d \tau.$
  Set $v(x,t)=u(x,t)+\mathcal{F}_2(t)$. Then
    $$\widetilde L v(x,t)=\widetilde  Lu(x,t) + \mathcal{F}_2'(t)=Lu(x,t)+ F_2(t)\ge 0.$$
   On $\Gamma_t$, one has $v(x,t)\ge u(x,t)\ge G_2(t)$. By the Maximum Principle for $\widetilde L$ and $v$, we have  $v(x,t)\ge G_2(t)$ on $U_T$. Thus, $u(x,t)\ge G_2(t)-\mathcal{F}_2(t)$ on $U_T$,  which proves \eqref{min2}. 

\medskip
 \ref{MP3} Define, for $t\in[0,T]$, $\mathcal{F}(t)=\int_0^t F(t)\d \tau.$ 
 Applying (i) and (ii) to $F_1=F$ and $F_2=F$, respectively, we have
 \beqs
 u(x,t)\le G_1(t)+\mathcal{F}(t)\le G(t)+\mathcal{F}(t)\text{ and } 
 u(x,t)\ge G_2(t)-\mathcal{F}(t)\ge -G(t)-\mathcal{F}(t),
 \eeqs
 thus, obtain \eqref{maxabs}.

\medskip
We consider Case \ref{mb} now. 
Set $\widehat L w= w_t -\langle A,D^2 w\rangle + \widetilde b\cdot \nabla w$, where
    $$\widetilde b(x,t)= u(x,t) B(x,t)+ P'(u(x,t)) K(x,t)\nabla u(x,t).$$
We have $Lu=\widehat L u$.
By applying the results for Case \ref{ma} to $\mathscr L:=\widehat L$, we obtain \ref{MP1}, \ref{MP2} and \ref{MP3} again.
\end{proof}

The following is a simple consequence of Theorem \ref{maxprin2}. It shows that we can find preliminary estimates solutions of \eqref{nonIBVP} for all time  in terms of $f$, $g$ and $u_0$. These simple estimates will be improved much more later. Actually, they  play the role of a starting point for those improvements. 
 
\begin{corollary}\label{range}
Let $A(x,t)$ satisfy \eqref{Aelip}  with $T=\infty$ and let $u\in C(\overline{U_\infty})\cap C_{x,t}^{2,1}(U_\infty)$ be a solution of \eqref{nonIBVP}. Assume $f^+(x,t)\le F_1(t)$, resp., $f^-(x,t)\le F_2(t)$ for all $(x,t)\in U_\infty$, where $F_1$, resp., $F_2$  is an integrable function in $C([0,\infty),[0,\infty))$,
and $g(x,t)$ is bounded from above, resp., from below on $\Gamma\times (0,\infty)$.
  Denote 
    \beq \label{mtm1}
  m_1= \max \left \{\sup_{x\in U} u_0(x),\sup_{(x,t)\in \Gamma\times (0,\infty)} g(x,t)\right\}
  \text{ and } \widetilde  m_1= \int_0^\infty F_1(\tau)\d\tau,
    \eeq
resp.,
  \beq \label{mtm2}
  m_2= \min \left \{\inf_{x\in U} u_0(x),\inf_{(x,t)\in \Gamma\times (0,\infty)} g(x,t)\right\}
  \text{ and } \widetilde  m_2=   \int_0^\infty F_2(\tau)\d\tau.
    \eeq
 Then one has 
    \beq\label{uuper}
    u(\overline{U_\infty})\subset (-\infty,m_1+\widetilde  m_1], \text{ resp., }
    u(\overline{U_\infty})\subset [m_2- \widetilde m_2,\infty).
    \eeq
\end{corollary}
\begin{proof}
For $(x,t)\in U_\infty$, we have 
\beqs 
Lu(x,t)=f(x,t)\le f^+(x,t)\le F_1(t),
\text{ resp., }
Lu(x,t)=f(x,t)\ge -f^-(x,t)\ge -F_2(t).
\eeqs 
For $t>0$, applying inequality \eqref{max2}, resp., \eqref{min2} in Theorem \ref{maxprin2} to $T:=2t$, we obtain
\beqs
        \max_{\overline{U_t}} u \le \max_{\Gamma_t} u+\int_0^t F_1(\tau)\d \tau
         \le \sup_{\Gamma_\infty} u+\int_0^\infty F_1(\tau)\d \tau
        = m_1 + \widetilde m_1,
\eeqs
    resp.,
\beqs
        \min_{\overline{U_t}} u \ge \min_{\Gamma_t} u-\int_0^t F_2(\tau)\d \tau
         \ge \inf_{\Gamma_\infty} u-\int_0^\infty F_2(\tau)\d \tau
        = m_2 -\widetilde m_2.
\eeqs
    Thus, the statement \eqref{uuper} follows.
\end{proof}

\section{The linear problem}\label{linsec}
Let $T$ be an extended number in $(0,\infty]$, and 
\beq \label{Abgive}
A:U_T\to \mathcal M^{n\times n}_{{\rm sym}}
\text{ and }b:U_T\to \R^n \text{ be given functions.}
\eeq 
Define the linear operator $\widetilde L $ by \eqref{Ltil}.

\begin{assumption}\label{firstA}
    There are constants $c_0>0$ and $M_1>0$ such that \eqref{Aelip} holds and 
\beqs
{\rm Tr}(A(x,t))\le M_1 \text{ for all }(x,t)\in U_T.
\eeqs
\end{assumption}

\begin{assumption}\label{condB}
There is a constant $M_2\ge 0$ so that
\beqs
 |b(x,t)|\le M_2 \text{ for all }(x,t)\in U_T.
\eeqs
\end{assumption}

We recall a Growth Lemma of Landis-type \cite{LandisBook}.

\begin{lemma}[Growth Lemma {\cite[Lemma 4.2]{HI3}}] \label{lemgrowth}
Let Assumptions \ref{firstA} and \ref{condB} hold for some number $T\in(0,\infty)$. Given a point $x_0\not\in \bar U$, let $r_0=r_*(x_0)$ and $R=R_*(x_0)$. Set 
 \beq\label{sTe}
 \beta_*=\frac{M_1+M_2R}{2c_0},\  
 \beta=\max\left \{ \beta_*,\frac{R^2}{4c_0T}\right\}, \  
 T_*=\frac{R^2}{4c_0 \beta}\in(0,T],\   
 \eta_*=1-(r_0/R)^{2\beta}\in(0,1).
 \eeq
Let 
\beq\label{wclass} w\in C(\overline{U_T})\cap C_{x,t}^{2,1}(U_T).
\eeq
If $\widetilde Lw\le 0$ on $U_T$ and $w\le 0$ on $\Gamma\times[0,T]$, then one has
\beqs \label{Tgrow}
\max_{x\in\bar U} w^+(x,T_*)
\le \eta_* \max_{x\in\bar U} w^+(x,0).
\eeqs 
\end{lemma}

Notice in \eqref{sTe} that if 
\beq\label{Tinval} 
4c_0\beta_*T= R^2,
\eeq 
then 
\beq \label{easysTe}
\beta=\beta_*=R^2/(4c_0T),\quad 
T_*=R^2/(4c_0\beta_*)=T\text{ and } 
\eta_*=1-(r_0/R)^{2\beta_*}.\eeq 

\begin{lemma}[Growth Lemma, inhomogeneous version] \label{lemG2}
We assume the same as in Lemma \ref{lemgrowth} up to \eqref{wclass} and also assume \eqref{Tinval} holds.
Suppose  there is a function $F\in C([0,T],[0,\infty))$ such that 
    $\widetilde Lw(x,t)\le F(t)$ for all $(x,t)\in U_T$.     
Then one has
 \beq \label{Tgrow2}
\max_{x\in\bar U} w^+(x,T)
\le \eta_* \max_{x\in\bar U} w^+(x,0)
+\max_{(x,t)\in \Gamma\times [0,T]} w^+(x,t)+\int_{0}^{T}F(t)\d t.
\eeq 
\end{lemma}
\begin{proof}
Denote $ D=\int_{0}^{T}F(t)\d t,$ 
$ \delta=\max_{\Gamma\times [0,T]} w^+(x,t)$,   
$\Lambda=D+\delta$, and 
\beqs 
J(t)=\max_{x\in\bar U} w^+(x,t) \text{ for $t\in[0,T]$.}  
\eeqs 
Define the function 
\beqs 
v(x,t)=w(x,t)-\delta-\int_{0}^t F(\tau)\d\tau \text{ for }(x,t)\in \bar U\times[0,T].
\eeqs 
Then one has $v\in C(\overline{U_T})\cap C_{x,t}^{2,1}(U_T)$ and  
\beq \label{GLv}
\widetilde L v=\widetilde Lw -F(t)\le 0\text{ in }U_T.
\eeq 
Moreover, for the boundary values of $v$, one has
\beqs
v(x,t)\le w(x,t)-\delta\le w^+(x,t)-\delta\le 0
\text{ on }\Gamma\times[0,T].
\eeqs 
Note that we now have the values of $\beta$, $T_*$ and $\eta_*$ in \eqref{easysTe}.
By the previous Growth Lemma - Lemma \ref{lemgrowth} - applied to the solution $v$ of \eqref{GLv}, we have 
\beq\label{Gvmax1}
\max_{x\in\bar U} v^+(x,T)\le \eta_* \max_{x\in\bar U} v^+(x,0).
\eeq
On the one hand, $v(x,0)= w(x,0)-\delta\le w(x,0)$, 
hence, 
\beq\label{GvJ}
\max_{x\in\bar U} v^+(x,0)\le \max_{x\in\bar U} w^+(x,0)= J(0).
\eeq 
On the other hand,  $v(x,T)=w(x,T)-\delta-D=w(x,T)-\Lambda$, hence,
\beqs 
w(x,T)=v(x,T)+\Lambda\le \max_{x\in\bar U} v^+(x,T) +\Lambda.
\eeqs
Thus,
\beq\label{GJb}
J(T)\le  \max_{x\in\bar U} v^+(x,T) +\Lambda.
\eeq
Combining \eqref{GJb} with \eqref{Gvmax1} and \eqref{GvJ},  we obtain $J(T)\le \eta_* J(0)+\Lambda$ 
which proves \eqref{Tgrow2}.
\end{proof}

We next apply  Lemma \ref{lemG2} to consecutive intervals to obtain estimates, though still in very technical forms,  for all time.

\begin{proposition}\label{genlem}
Let $T=\infty$ and assume the function $A(x,t)$ satisfies Assumption \ref{firstA}. 
Given a strictly increasing sequence $(T_k)_{k=0}^\infty$ of nonnegative numbers, let $\tau_k=T_k-T_{k-1}$ for $k\ge 1$.
For each $k\ge 1$, supposing  
\beq \label{bmk}
|b(x,t)|\le m_k\text{ on }U\times (T_{k-1},T_k]
\eeq 
and letting $x_k$ be any point in $\R^n\setminus \bar U$, we set   
    \beq\label{rRbe}
    r_k=r_*(x_k),\  R_k=R_*(x_k),\ 
    \beta_k=\frac{M_1+m_k R_k}{2c_0}\text{ and } \eta_k=1-(r_k/R_k)^{2\beta_k}.
    \eeq 
Assume 
    \beq\label{btkcond}
    4c_0\beta_k \tau_k= R_k^2 \text{ for all }k\ge 1.
    \eeq 

Let $w\in C(\bar U\times[T_0,\infty))\cap C_{x,t}^{2,1}(U\times(T_0,\infty))$ and  define 
\beq\label{Jkdef}
J_k=\max_{x\in \bar U} w^+(x,T_k)\text{ for $k\ge 0$.}
\eeq  
Suppose there is a function $F\in C([T_0,\infty),[0,\infty))$ such that 
  \beq\label{LwF}   
  \widetilde Lw(x,t)\le F(t)\text{ for all } (x,t)\in U\times(T_0,\infty).
  \eeq 
    
Then one has, for all $k\ge 1$,
\beq\label{wv1}
\max_{(x,t)\in \bar U\times[T_{k-1},T_k]}w^+(x,t)
\le J_{k-1}+\Lambda_k
\eeq
and
\beq\label{Jkest0}
J_k\le \sum_{m=0}^k \left[ \left(\prod_{j=m+1}^k \eta_j\right) \Lambda_m\right],
\eeq
where
\beq \label{lamk}
\Lambda_0=J_0\text{ and }
\Lambda_k=\max_{(x,t)\in \Gamma\times [T_{k-1},T_k]} w^+(x,t)+\int_{T_{k-1}}^{T_k}F(t)\d t
 \text{ for $k\ge 1$.  }
\eeq 
\end{proposition}
\begin{proof}
Thanks to \eqref{LwF}, we can apply, for each $k\ge 1$, Theorem \ref{maxprin2}, case \ref{ma}, part \ref{MP1} to the interval $[0,\tau_k]$ in place of $[0,T]$, functions $u(t):=w(T_{k-1}+t)$, $F_1(t):=F(T_{k-1}+t)$ and select $t=\tau_k$ in \eqref{max2}. It results in        
\begin{align*}
        \max_{\bar U\times [T_{k-1},T_k]} w 
        &\le \max_{(\bar U\times\{T_{k-1}\})\cup (\Gamma \times [T_{k-1},T_k])} w 
        +\int_{T_{k-1}}^{T_k} F(\tau)\d \tau \\
        &\le \left (J_{k-1}+\max_{\Gamma \times [T_{k-1},T_k]} w^+ \right)+\int_{T_{k-1}}^{T_k} F(\tau)\d \tau
        \le J_{k-1}+\Lambda_k
\end{align*}
which implies \eqref{wv1}.

Let $k\ge 1$. We apply Lemma \ref{lemG2} to the interval $[T_{k-1},T_k]$ in place of $[0,T]$, with the same $M_1$, but, referring to the numbers in \eqref{bmk} and \eqref{rRbe}, 
\beqs 
M_2:=m_k,\  
x_0:=x_k,\ 
r_0:=r_k,\ R:=R_k,\ 
\beta=\beta_*:=\beta_k,\ 
\eta_*:=\eta_k.
\eeqs 
We obtain from \eqref{Tgrow2} that
\beqs
J_k\le \eta_k J_{k-1}+\Lambda_k.
\eeqs 
Iterating this estimate for $k-1$, $k-2$, etc., we have
\begin{align*}
J_k&\le \eta_k(\eta_{k-1}J_{k-2}+\Lambda_{k-1})+\Lambda_k
=\eta_k\eta_{k-1}J_{k-2}+\eta_k\Lambda_{k-1}+\Lambda_k\\
&\le \eta_k\eta_{k-1}\eta_{k-2}J_{k-3}+\eta_k\eta_{k-1}\Lambda_{k-2}+\eta_k\Lambda_{k-1}+\Lambda_k\\
&\le \ldots \le \eta_k\eta_{k-1}...\eta_1 J_0 + \sum_{m=1}^k  \prod_{j=m+1}^k \eta_j \Lambda_m,
\end{align*}
which implies \eqref{Jkest0}.
\end{proof}

Below, we give more specific descriptions of the behavior of the solution $w(x,t)$, for large $t$ in the case $T=\infty$, of the equation $\widetilde Lw(x,t)=f(x,t)$. 
For the rest of this section,  $A$ and $b$ are fixed functions in \eqref{Abgive} with $T=\infty$,  and  the linear operator $\widetilde L$ is correspondingly defined by \eqref{Ltil} with $T=\infty$.

\subsection{Case of bounded drifts}\label{bdedsec}

We first study the case when the drift term $b(x,t)$ is uniformly bounded in $x$ and $t$.
We will apply Proposition \ref{genlem} which raises the issue of estimating the series in \eqref{Jkest0}. In some particular cases, we will use the following simple result.

\begin{lemma}\label{dH}
    Given a number  $\eta\in(0,1)$ and a sequence $(\Lambda_k)_{k=0}^\infty$ of nonnegative numbers. For $k\ge 0$, let 
$ a_k= \sum_{j=0}^k \eta^{k-j}\Lambda_j$.
Then
\beq\label{alim}
\limsup_{k\to\infty} a_k\le \frac1{1-\eta} \limsup_{k\to\infty}\Lambda_k.
\eeq
\end{lemma}
\begin{proof}
If $(\Lambda_k)_{k=0}^\infty$ is unbounded, then inequality \eqref{alim} holds true. 
Now assume there is $\Lambda_*>0$ such that $\Lambda_k\le \Lambda_*$ for all $k\ge 0$.
For $m\ge 0$, let $\varep_m=\sup\{\Lambda_j:j\ge m\}$.
Then, for $k>m$, we have
\begin{align*}
    a_k&=\sum_{j=0}^m \eta^{k-j}\Lambda_j +\sum_{j=m+1}^k \eta^{k-j}\Lambda_j \le \Lambda_* \sum_{j=0}^m \eta^{k-j}+\varep_m \sum_{j=m+1}^k \eta^{k-j}\\
    &    \le  \Lambda_* \eta^k\left (\sum_{j=0}^m \eta^{-j}\right )+ \frac{\varep_m}{1-\eta}.
\end{align*}
Passing $k\to\infty$ yields $\displaystyle \limsup_{k\to\infty} a_k\le \varep_m/(1-\eta).$ 
Then taking $m\to \infty$, we obtain \eqref{alim}.
\end{proof}

Our first explicit asymptotic estimates are the following.

\begin{theorem}\label{NHthm}
Under  Assumptions \ref{firstA} and \ref{condB} both for $T=\infty$, let $w\in C(\overline{U_\infty})\cap C_{x,t}^{2,1}(U_\infty)$. 
Given any point $x_0\not\in \bar U$, set the numbers $r_0=r_*(x_0)$, $R=R_*(x_0)$,
\beq\label{beet}
 \beta_*=\frac1{2c_0}(M_1+M_2R), \
 T_*=\frac{R^2}{4c_0 \beta_*}\text{ and }
 \eta_*=1-(r_0/R)^{2\beta_*}.
\eeq

\begin{enumerate}[label=\tnum]
    \item \label{dec1}    Suppose there is a functions $F_1\in C([0,\infty),[0,\infty))$ such that 
    $\widetilde Lw(x,t)\le F_1(t)$ for all $(x,t)\in U_\infty$.  
Then 
\beq\label{lsup1}
\limsup_{t\to\infty} \left[\max_{x\in\bar U} w^+(x,t)\right]
\le \frac{2-\eta_*}{1-\eta_*}\left(
\limsup_{t\to\infty}\left[\max_{x\in\Gamma} w^+(x,t)\right]
+\limsup_{t\to\infty}\int_t^{t+T_*} F_1(\tau)\d\tau 
\right).
\eeq 

\item \label{dec3}     Suppose there is a function $F\in C([0,\infty),[0,\infty))$ such that 
    $|\widetilde Lw(x,t)|\le F(t)$ for all $(x,t)\in U_\infty$. 
Then 
\beq\label{lasb1}
\limsup_{t\to\infty} \left[\max_{x\in\bar U} |w(x,t)|\right]
\le \frac{2-\eta_*}{1-\eta_*}\left(
\limsup_{t\to\infty}\left[\max_{x\in\Gamma} |w(x,t)|\right]
+\limsup_{t\to\infty}\int_t^{t+T_*} F(\tau)\d\tau 
\right).
\eeq
\end{enumerate}
\end{theorem}
\begin{proof}

\medskip\noindent
\ref{dec1} For $k\ge 0$, let $T_k=kT_*$
and 
$J_k=\max_{x\in \bar U} w^+(x,T_k)$.
Let  $\Lambda_0=J_0$, and, for $k\ge 1$,  
\beq \label{dDk}
\delta_k=\max_{\Gamma\times [T_{k-1},T_k]} w^+,\quad 
D_k=\int_{T_{k-1}}^{T_k}F_1(t)\d t,\quad 
\Lambda_k=D_k+\delta_k.
\eeq 
We apply  Proposition \ref{genlem} to $x_k=x_0$, same $M_1$, and, referring to \eqref{bmk}, $m_k=M_2$.
Then $\tau_k=T_{k+1}-T_k=T_*$,  the numbers in \eqref{rRbe} are 
$r_k=r_0$, $R_k=R$, $\beta_k=\beta_*$,  $\eta_k=\eta_*$
and condition \eqref{btkcond} is met.
It follows from \eqref{Jkest0} and Lemma \ref{dH} that 
\beq\label{Jk1}
\limsup_{k\to\infty} J_k
\le \limsup_{k\to\infty} \left( \sum_{m=0}^{k} \eta_*^{k-m} \Lambda_m\right)
\le \frac1{1-\eta_*}\limsup_{k\to\infty}\Lambda_k.
\eeq
By \eqref{wv1} and \eqref{Jk1},
\beq\label{wLam1}
\limsup_{t\to\infty}\left( \max_{x\in\bar U} w^+(x,t)\right)\le \limsup_{k\to\infty} J_{k-1}+\limsup_{k\to\infty} \Lambda_k
\le \frac{2-\eta_*}{1-\eta_*}\limsup_{k\to\infty} \Lambda_k.
\eeq
One has from \eqref{dDk} that  
\beq\label{LdD}
\limsup_{k\to\infty} \Lambda_k\le 
\limsup_{k\to\infty} \delta_k
+\limsup_{k\to\infty} D_k 
\le \limsup_{t\to\infty}\left(\max_{x\in\Gamma} w^+(x,t)\right)
+\limsup_{t\to\infty}\int_t^{t+T_*} F_1(\tau)\d\tau.
\eeq
Thus, by combining \eqref{wLam1} with \eqref{LdD} we obtain the desired inequality \eqref{lsup1}.

\medskip\noindent
\ref{dec3} It suffices to consider the case when the right-hand side of \eqref{lasb1} is a finite number which we denote by $C_0$.
Applying part (i) to function $F_1:=F$, we obtain
\beq\label{get1}
\limsup_{t\to\infty} \left[\max_{x\in\bar U} w^+(x,t)\right]
\le \frac{2-\eta_*}{1-\eta_*}\left(
\limsup_{t\to\infty}\left[\max_{\Gamma} w^+(x,t)\right]
+\limsup_{t\to\infty}\int_t^{t+T_*} F(\tau)\d\tau 
\right)\le C_0.
\eeq
Note that $\widetilde Lw\ge -F(t)$, hence $\widetilde L(- w)\le F(t).$
Applying part (i) to functions $(-w)$ and $F_1(t):=F(t)$, we obtain, recalling that $(-w)^+=w^-$,
\beq\label{get2}
\limsup_{t\to\infty} \left[\max_{x\in\bar U} w^-(x,t)\right]
\le \frac{2-\eta_*}{1-\eta_*}\left(
\limsup_{t\to\infty}\left[\max_{x\in\Gamma} w^-(x,t)\right]
+\limsup_{t\to\infty}\int_t^{t+T_*} F(\tau)\d\tau 
\right)\le C_0.
\eeq
 Observe that
\beqs 
\max_{x\in \bar U}|w(x,t)|
=\max_{x\in \bar U}\left[\max\{ w^+(x,t),w^-(x,t)\}\right]\le  \max\left \{ \max_{x\in \bar U} w^+(x,t),\max_{x\in \bar U} w^-(x,t)\right\}.
\eeqs
Then one has
 \beq\label{absw}
\limsup_{t\to \infty}\left[\max_{x\in \bar U} |w(x,t)|\right]
  \le \max\left \{
 \limsup_{t\to \infty}\left[\max_{x\in \bar U} w^+(x,t)\right], 
 \limsup_{t\to \infty}\left[\max_{x\in \bar U} w^-(x,t)\right]
 \right\}.
 \eeq
Combining this inequality \eqref{absw} with the  estimates \eqref{get1} and \eqref{get2}, 
one obtains \eqref{lasb1}.
\end{proof}

Theorem \ref{NHthm} directly yields a result on the decay of $|w(x,t)|$.

\begin{corollary}\label{NHCor}
Under  Assumptions \ref{firstA} and \ref{condB} both for $T=\infty$, 
let $w\in C(\overline{U_\infty})\cap C_{x,t}^{2,1}(U_\infty)$ and 
assume there is a function $F\in C([0,\infty),[0,\infty))$ such that 
    $|\widetilde Lw(x,t)|\le F(t)$ for all $(x,t)\in U_\infty$. 
If, in addition, 
\beq\label{nhwdata}
\lim_{t\to\infty}\left(\max_{x\in\Gamma} |w(x,t)|\right)=0\text{ and } 
\lim_{t\to\infty} \int_t^{t+1}F(\tau)\d \tau=0,
\eeq
    then 
  \beq   \label{nhwlim}
\lim_{t\to\infty} \left( \max_{x\in\bar U}|w(x,t)|\right)=0.
  \eeq
\end{corollary}
\begin{proof}
Let $T_*$ and $\eta_*$ be defined as in \eqref{beet}.
We fix a positive integer $N\ge T_*$.
Observe that
\beq\label{limso}
\begin{aligned}
\limsup_{t\to\infty}\int_t^{t+T_*}F(\tau)\d \tau 
& \le \limsup_{t\to\infty}\int_t^{t+N} F(\tau)\d \tau
    =\limsup_{t\to\infty}\sum_{j=0}^{N-1}\int_{t+j}^{t+j+1} F(\tau)\d \tau\\
&    \le N\limsup_{t\to\infty} \int_t^{t+1} F(\tau)\d \tau=0.
\end{aligned}
\eeq 
Then \eqref{nhwlim} follows from \eqref{lasb1}, \eqref{nhwdata} and \eqref{limso}.
\end{proof}

Note that the second condition in \eqref{nhwdata} involves only the integral of $F$ on the time interval $[t,t+1]$. It does not require $F(t)\to 0$ as $t\to\infty$.

\subsection{Case of unbounded drifts}\label{ubsec}

Below, we study the case when $\sup_{x\in U}|b(x,t)|$ is unbounded as $t\to\infty$. Let  $w\in C(\overline{U_\infty})\cap C_{x,t}^{2,1}(U_\infty)$.
The following technical conditions will be needed in our later proofs.

\begin{assumption}\label{abZL}
Assume there is $T^*\ge 0$ such that one has the following.

\begin{enumerate}[label=\tnum]
\item\label{FbL1} There is a function $F\in C([T^*,\infty),[0,\infty))$ such that 
    \beq \label{Lwf}
    \widetilde Lw(x,t)\le F(t)\text{ for all }(x,t)\in U\times(T^*,\infty).
    \eeq 

    \item \label{FbL2} There is a continuous function $z_0(t)>0$ on $[T^*,\infty)$ increasing to infinity as $t\to\infty$ such that 
\beq\label{bz}
|b(x,t)|\le z_0(t) \text{ for all }(x,t)\in U\times (T^*,\infty).
\eeq
    \item \label{FbL3} There is a decreasing, continuous function $\Lambda_*(t)\ge 0$ on $[T^*,\infty)$ such that
\beq\label{lamst}
\max_{(x,t)\in \Gamma\times[t,t+1]}w^+(x,t)
+\int_t^{t+1}F(\tau)\d \tau \le  \Lambda_*(t) \text{ for all }t\ge T^*.
\eeq 
\end{enumerate}
\end{assumption}

Referring to \eqref{lamst}, we define a function $\bar\Lambda(t)$ on $[4T^*,\infty)$ by 
\beq\label{lamdef}
\bar \Lambda(t)=\Lambda_*(t/4) \text{ for }t\ge 4T^*.
\eeq

\begin{assumption}\label{AssumZ}
Under Assumption \ref{abZL}, referring to the functions $z_0(t)$ in \eqref{bz}  and $\bar \Lambda(t)$ in \eqref{lamdef}, we assume
\beq\label{z2cond}
\int_{T^*}^{\infty} e^{-\frac{d}{c_0}z_0(x)} \d x =\infty,
\eeq
and  there are numbers $t_*\ge 4T^*$ and $\varep_*>0$ such that, for any $\varep\in(0,\varep_*)$, the function 
\beq \label{monocond}
t\in [t_*,\infty)\mapsto \bar \Lambda(t)\exp\left(e^{-\frac{d^2}{2c_0}-\varep}\int_{t_*}^{t+1} e^{-\frac{d}{c_0}z_0(x)} \d x\right)\text{  is  monotone.}
\eeq 
\end{assumption}

For the upper bound of $w(x,t)$, we have the next theorem.

\begin{theorem}\label{NHthm3}
Under Assumption \ref{firstA} with $T=\infty$ and Assumptions \ref{abZL} and \ref{AssumZ}, if
\beq\label{Llim}
\limsup_{t\to\infty} \left(\bar \Lambda(t) e^{\frac{d}{c_0}z_0(t)}\right)=\ell\in [0,\infty),
\eeq
then $w(x,t)$ is bounded from above on $\overline{U_\infty}$ and 
\beq\label{lel1}
\limsup_{t\to\infty}\left( \max_{x\in\bar U} w^+(x,t)\right)\le e^\frac{d^2}{2c_0} \ell .
\eeq    
\end{theorem}

\begin{remark} \label{shiftime} The following remarks are in order.
\begin{enumerate}[label=\rnum]
    \item     If Assumption \ref{abZL} holds for some $T^*\ge 0$, then, for any $\bar T>T^*$, Assumption \ref{abZL} also holds with $T^*$ being replaced with $\bar T$.
    
    \item If Assumption \ref{AssumZ} holds for some $T^*\ge 0$ and $t_*\ge 4T^*$, then, for any $\bar T\ge T^*$ and $\bar t\ge 4\bar T$, Assumption \ref{AssumZ} also holds with $T^*$ and $t_*$  being replaced by $\bar T$ and $\bar t$ respectively.

    \item The limit superior \eqref{Llim} does not depend on the choice of $T^*$ and $t_*$ in Assumptions \ref{abZL} and \ref{AssumZ}.
    \end{enumerate}
\end{remark}

The proof of Theorem \ref{NHthm3} involves constructing barrier functions with different centers.
These centers will be selected to, at least, optimize the inequalities in \eqref{rRd}.
In fact, we have the following simple but useful fact.

\begin{lemma}\label{choicey}
    For any $R>d$, there exists a point $y\not \in \bar U$ such that
    \beq\label{rexact}
    r_*(y)=R-d\text{ and }R_*(y)=R.
    \eeq
\end{lemma}
\begin{proof}
Let  $\xi_1,\xi_2\in \bar U$ such that $|\xi_1-\xi_2|=d>0$. 
First, let $y=y(t)\eqdef \xi_1 +t(\xi_2-\xi_1)$ with number $t> 1$. 
We have  $|y-\xi_1|=td$ and $|y-\xi_2|=(t-1)d$, which yield $|y-\xi_1|=d+|y-\xi_2|$.

Because $t>1$, one has $|y-\xi_1|>d$ which implies $y\not\in \bar U$.

For $z\in \bar U$, $|y-z|\le |y-\xi_2|+|\xi_2-z|\le |y-\xi_2|+d=|y-\xi_1|$. Hence, $|y-\xi_1|=R_*(y)$.

For $z\in \bar U$, $|y-z|\ge |y-\xi_1|-|\xi_1-z|\ge |y-\xi_1|-d=|y-\xi_2|$. Hence, $|y-\xi_2|=r_*(y)$.

Then the point $y=y(R/d)$ is not in $\bar U$ and it satisfies \eqref{rexact}.
\end{proof}

\begin{proof}[Proof of Theorem \ref{NHthm3}]
Because of the continuity of $w$ on $\overline{U_\infty}$, it is bounded on $\overline{U_T}$ for any number $T\in(0,\infty)$.
Thus, the estimate \eqref{lel1}, once established, will imply that $w(x,t)$ is bounded from above on $\overline{U_\infty}$.
Therefore, it suffices to prove \eqref{lel1} which we proceed in many steps below. 

\medskip\noindent
\textbf{Step 1.}  Define the function 
\beq\label{Zdef}
Z(t)= -\frac{M_1+z_0(t)^2}{c_0}\ln \left(1-\frac{d}{z_0(t)}\right) \text{ for }t\ge T^*.
\eeq
By the Taylor expansion, we have
\beqs
-\ln(1-x)=x+\frac12x^2+\mathcal O(x^3)\text{ as } x\to 0.
\eeqs
Consequently,  as $t\to\infty$,
\begin{align}
Z(t)
&=\frac{z_0(t)^2}{c_0} \left( \frac{d}{z_0(t)}+\frac12\left(\frac{d}{z_0(t)}\right)^2+\mathcal O(1/z_0(t)^3)\right)
-\frac{M_1}{c_0}\ln \left(1-\frac{d}{z_0(t)}\right) \notag\\
&=\frac{d}{c_0}z_0(t)+\frac{d^2}{2c_0}+\mathcal O(1/z_0(t))+\mathcal O(1/z_0(t)).\label{Zex}
\end{align}

Let $\varep_0$ be any number in the interval $(0,\varep_*)$.
Thanks to \eqref{Zex}, we can take a sufficiently large number $T_0>\max\{1/3,t_*\}$ such that 
\beq\label{T0choice}
z_0(T_0)>\max\left\{d,\sqrt{M_1}\right\}  \text{ and }
Z(t)\le \frac{d}{c_0}z_0(t)+\frac{d^2}{2c_0}+\varep_0\text{ for all }t\ge T_0.
\eeq
For $t\ge 0$, define  
\beq\label{zz} 
z_*(t)=z_0(T_0+t)\text{ and } 
z(t)= \frac{d}{c_0}z_*(t)+\frac{d^2}{2c_0}+\varep_0
= \frac{d}{c_0}z_0(T_0+t)+\frac{d^2}{2c_0}+\varep_0 .
\eeq
We then rewrite the inequality for $Z(t)$ in \eqref{T0choice} as 
\beq\label{Zprop}
Z(T_0+t)\le 
\frac{d}{c_0}z_0(T_0+t)+\frac{d^2}{2c_0}+\varep_0
=z(t)\text{ for all }t\ge 0.
\eeq
Note that the functions $z_0(t)$ and $z(t)$ are increasing and, for $t\ge 0$, 
\beq \label{zd}
z_*(t)\ge z_0(T_0)>d.
\eeq 
\medskip\noindent
\textbf{Step 2.} For each integer $k\ge 1$, by \eqref{zd} and the virtue of  Lemma \ref{choicey}, we select a point $x_k\not \in \bar U$ so that
\beq\label{rRk}
R_k :=R_*(x_k)=z_*(k), \quad r_k:=r_*(x_k)=R_k-d,
\eeq
and define 
\beqs
\widetilde\beta_k=\frac{M_1+R_k z_*(k)}{2c_0},\quad 
\tau_k=\frac{R_k^2}{4c_0 \widetilde\beta_k}=\frac{R_k^2}{2(M_1+R_k z_*(k))}.
\eeqs
For $k\ge 1$, set $T_k=T_0+\sum_{j=1}^k\tau_j$. 
We note from the choice \eqref{rRk} that
\beq\label{tausmall}
\tau_k\le  \frac{R_k^2}{2R_kz_*(k)}=\frac12.
\eeq
Also, by \eqref{T0choice}, the increase of $z_*(t)$ and \eqref{rRk}, $M_1\le z_0(T_0)^2=z_*(0)^2\le z_*(k)^2=R_k^2$, hence
\beqs
\tau_k=\frac{R_k^2}{2(M_1+R_k z_*(k))}\ge \frac{R_k^2}{2(R_k^2+R_k^2)}=\frac14.
\eeqs
Therefore, for all $k\ge 0$,
\beq\label{Tran}
T_0+k/4\le T_k\le T_0+k/2.
\eeq

\medskip\noindent
\textbf{Step 3.} 
We will apply Proposition \ref{genlem}.
Let $k\ge 1$. Noticing from \eqref{Tran} that $[T_{k-1},T_k]\subset [T_0,T_0+k]$, hence,  we have  on $U\times (T_{k-1},T_k]$
\beqs
|b(x,t)|\le z_0(t)\le z_0(T_k)\le z_0(T_0+k)=z_*(k).
\eeqs
This prompts us to  define $m_k=z_*(k)$ so that inequality \eqref{bmk} holds.
Referring to  definitions of the numbers $\beta_k$ and $\eta_k$ in Proposition \ref{genlem}, we have
\beq\label{bebek}
\beta_k=\frac{M_1+R_k m_k}{2c_0}=\frac{M_1+R_k z_*(k)}{2c_0}=\widetilde \beta_k,
\eeq
\beq \label{etak}
\eta_k=1-(r_k/R_k)^{2\beta_k}=1-(r_k/R_k)^{2\widetilde\beta_k}.
\eeq
Clearly, $4c_0 \beta_k \tau_k=4c_0\widetilde \beta_k \tau_k=R_k^2$,
hence condition \eqref{btkcond} is met.

\medskip\noindent
\textbf{Step 4.} For $k\ge 1$, we write $\eta_k$ in \eqref{etak}, with the use of \eqref{rRk} and \eqref{bebek}, explicitly as
\beqs
\eta_k
=1-\left(\frac{z_*(k)-d}{z_*(k)}\right)^\frac{M_1+z_*^2(k)}{c_0}
= 1-e^{-\tilde z(k)},
\eeqs
where
\beqs
\tilde z(t)= -\frac{M_1+z_*(t)^2}{c_0}\ln \left(1-\frac{d}{z_*(t)}\right) \text{ for $t\ge 0$.}
\eeqs
Note from \eqref{Zdef} and \eqref{zz} that $\tilde z(t)=Z(T_0+t)$. 
This identity and  and property \eqref{Zprop} imply $\tilde z(t)\le z(t)$ for all $t\ge 0$.
Therefore, 
\beq\label{etexz}
\eta_k\le  1-e^{-z(k)} \text { for $k\ge 1$.}
\eeq

\medskip\noindent
\textbf{Step 5.}
Let $J_k$ and $\Lambda_k$ be defined by \eqref{Jkdef} and \eqref{lamk} respectively.
For $k\ge 1$, by using the fact $\tau_k\le 1$ in \eqref{tausmall}, we have
\beq \label{lamk1}
\Lambda_k\le 
\max_{\Gamma\times [T_{k-1},T_{k-1}+1]} w^+(x,t)
+\int_{T_{k-1}}^{T_{k-1}+1}F(t)\d t
\le \Lambda_*(T_{k-1}).
\eeq 
For $k\ge 1$, estimating $T_{k-1}$ by the first inequality in \eqref{Tran}, and using the facts that $\Lambda_*(t)$ is decreasing and $3T_0-1>0$, we have 
\beq \label{Lamineq}
\Lambda_*(T_{k-1})\le \Lambda_*\left(T_0+\frac{k-1}4\right)=\Lambda_*\left(\frac{T_0+k+3T_0-1}4\right)
\le \Lambda_*\left(\frac{T_0+k}4\right).    
\eeq 
Define the function $\Lambda(t)$ by 
\beqs 
\Lambda(t)=\bar \Lambda(T_0+t)
=\Lambda_*\left(\frac{T_0+t}4\right)\text{ for $t\ge 0$.}
\eeqs 
Then \eqref{lamk1} and \eqref{Lamineq} yield 
\beq\label{lamlam}
\Lambda_k\le \Lambda(k) \text{ for all }k\ge 1.
\eeq
Thanks to \eqref{Llim} and the fact $z_0(t)\to\infty$ as $t\to\infty$, one has $\Lambda(t)\to 0$ as $t\to\infty$. As a consequence of this and \eqref{lamlam}, one has
\beq\label{lamzero}
\lim_{k\to\infty} \Lambda_k =0.
\eeq
\medskip\noindent
\textbf{Step 6.} 
We now apply Proposition \ref{genlem}.
By \eqref{wv1} and \eqref{lamzero},
\beq\label{wJk}
\limsup_{t\to\infty}\left( \max_{x\in\bar U} w^+(x,t)\right)\le \limsup_{k\to\infty} J_{k-1}+\limsup_{k\to\infty} \Lambda_k
\le \limsup_{k\to\infty} J_{k}.
\eeq
By inequality \eqref{Jkest0}, we have, for $k\ge 1$,
\beq\label{reJk}
J_k\le \sum_{m=0}^k  (I_{k,m}\Lambda_m),\text{ where } I_{k,m}=\prod_{j=m+1}^k \eta_j.
\eeq

Let $k\ge 1$ and $0\le m\le k$. Using \eqref{etexz} and inequality $\ln (1+x)\le x$,  for $x\in(-1,\infty)$, we have
\begin{align*}
\ln I_{k,m}
\le \sum_{j=m+1}^k \ln \left(1-e^{-z(j)}\right)
\le -\sum_{j=m+1}^k e^{-z(j)}.
\end{align*}
Because $z(j)$ is increasing in $j$, so is $-e^{-z(j)}$ and, hence, we can estimate 
$$-e^{-z(j)}\le -\int_{j}^{j+1} e^{-z(x)} \d x.$$
Thus, we obtain
\begin{align*}
\ln I_{k,m}
&
\le -\sum_{j=m+1}^k \int_{j}^{j+1} e^{-z(x)} \d x
= -\int_{m+1}^{k+1} e^{-z(x)} \d x.
\end{align*}
Setting
$E_0(t)=\exp\left(\int_{0}^{t} e^{-z(x)} \d x\right) \text{ for }t\ge 0$,
we have $\ln I_{k,m}\le E_0(m+1)-E_0(k+1)$. Therefore, 
\beqs
I_{k,m} \le \frac{E_0(m+1)}{E_0(k+1)}, \text{ hence, }
\sum_{m=0}^k I_{k,m} \Lambda_m
\le \frac{1}{E_0(k+1)}
\sum_{m=0}^k \left[ \Lambda_m E_0(m+1)\right].
\eeqs 
In the last sum, we consider  $m=0$ separately and then use the relation \eqref{lamlam} for $m\ge 1$.
It yields that
\beq\label{spl1}
\begin{aligned}
\sum_{m=0}^k I_{k,m} \Lambda_m
&\le \frac{\Lambda_0E_0(1)}{E_0(k+1)}
+\frac{1}{E_0(k+1)}
\sum_{m=1}^k \left[ \Lambda(m)E_0(m+1)\right].    
\end{aligned}
\eeq
By \eqref{zz},
\beqs
E_0(t)=\exp\left(\int_{0}^{t} e^{-\frac{d}{c_0}z_0(T_0+x)}e^{-\frac{d^2}{2c_0}-\varep_0 } \d x\right) 
=\exp\left( e^{-\frac{d^2}{2c_0}-\varep_0 } \int_{T_0}^{T_0+t} e^{-\frac{d}{c_0}z_0(x)} \d x\right).
\eeqs
Thanks to \eqref{z2cond}, we have $E_0(t)\to\infty$ as $t\to\infty$. 
Therefore, the first term on the right-hand side of \eqref{spl1} goes to zero as $k\to \infty$.
Combining this with \eqref{reJk}, we have
\beq\label{liJsum}
\limsup_{k\to\infty}J_k \le \limsup_{k\to\infty}\sum_{m=0}^k  I_{k,m}\Lambda_m
\le \limsup_{k\to\infty} \widetilde J_k,
\eeq 
where
\beq\label{tJ}
\widetilde J_k=\frac{1}{E_0(k+1)}
\sum_{m=1}^k \left[ \Lambda(m)E_0(m+1)\right].
\eeq
Observe, for $t\ge 0$,  that
\begin{align*} 
 \Lambda(t)E_0(t+1)
&=\Lambda(t)\exp \left(e^{-\frac{d^2}{2c_0}-\varep_0}\int_{0}^{t+1} e^{-\frac{d}{c_0} z_0(T_0+x)} \d x\right) \\
&=\bar \Lambda(T_0+t) \exp \left( e^{-\frac{d^2}{2c_0}-\varep_0}\int_{T_0}^{T_0+t+1} e^{-\frac{d}{c_0}z_0(x)} \d x\right).
\end{align*}
Hence,
 \beq  \label{iden}
\Lambda(t)E_0(t+1)= \frac{ \bar \Lambda(T_0+t) \exp \left( e^{-\frac{d^2}{2c_0}-\varep_0}\int_{t_*}^{T_0+t+1} e^{-\frac{d}{c_0}z_0(x)} \d x\right) }{\exp \left(e^{-\frac{d^2}{2c_0}-\varep_0}\int_{t_*}^{T_0} e^{-\frac{d}{c_0}z_0(x)} \d x\right)}.
\eeq 

\medskip\noindent
\textbf{Step 7.}
Regarding condition \eqref{monocond}, we first consider the case  the function 
\beq \label{increase}
q(t):=\bar \Lambda(t) \exp\left(e^{-\frac{d^2}{2c_0}-\varep_0}\int_{t_*}^{t+1} e^{-\frac{d}{c_0} z_0(x)} \d x\right)\text{ is increasing on $[t_*,\infty)$.}
\eeq 
Denote the positive constant in the denominator in \eqref{iden} by $E_*$.
We rewrite \eqref{iden} as
\beq\label{newiden}
\Lambda(t)E_0(t+1)= \frac{q(T_0+t)}{E_*} \text{ for }t\ge 0.
\eeq 
Therefore,  $\Lambda(t)E(t+1)$ is increasing. We use  this fact to bound the summand in \eqref{tJ} from above by
\beq\label{uplamE}
\Lambda(m)E_0(m+1)\le \int_{m}^{m+1} E_0(s+1) \Lambda(s) \d s.
\eeq 
From \eqref{tJ} and \eqref{uplamE}, we obtain
\beq \label{tJe}
\begin{aligned}
\widetilde J_k
&\le\frac{ \sum_{m=1}^k \int_{m}^{m+1} E_0(s+1) \Lambda(s) \d s }{E_0(k+1)}=\frac{\int_{1}^{k+1} E_0(s+1) \Lambda(s) \d s }{E_0(k+1)}\\
&
= \int_{1}^{k+1} \exp\left(-\int_{s+1}^{k+1} e^{-z(x)} \d x\right) \Lambda(s) \d s 
\le y(k+1),
\end{aligned}
\eeq 
where
\beq\label{ydef}
y(t)=\int_{0}^{t} \exp\left(-\int_{s+1}^t e^{-z(x)} \d x\right) \Lambda(s) \d s .
\eeq
Rewrite
\beqs
y(t)=\int_{0}^{t} \exp\left(-\int_{s}^{t} e^{-z(x)} \d x\right) f(s) \d s,\text{ where }
f(s)=\exp\left(\int_{s}^{s+1} e^{-z(x)} \d x\right) \Lambda(s).
\eeqs
We have, for $t>0$,
\beq\label{yeq}
y'(t)=-e^{- z(t)}y(t)+f(t)=-h(t)\varphi^{-1}(y(t))+f(t),
\eeq 
where $\varphi(x)=x$ and $h(t)=e^{- z(t)}$.
Note from condition \eqref{z2cond} that
$\int_0^\infty h(\tau)\d\tau =\infty.$
By applying \cite[Lemma A.1]{HIKS1} to equation \eqref{yeq}, we obtain 
\begin{align*}
\limsup_{t\to\infty}y(t)
&\le \varphi\left (\limsup_{t\to\infty }\left[ f(t)/h(t)\right]\right)
= \limsup_{t\to\infty }\left[ f(t) e^{ z(t)}\right]\\
&=\limsup_{t\to\infty }\left[ \Lambda(t) e^{z(t)}\exp\left(\int_{t}^{t+1} e^{-z(x)} \d x\right)\right].    
\end{align*}
Since $z(x)\to\infty $ as $x\to\infty$, one has
$\int_{t}^{t+1} e^{-z(x)} \d x\to 0$ as $t\to \infty$.
Consequently,
\beq\label{limy}
\limsup_{t\to\infty}y(t)
\le \limsup_{t\to\infty} \left( \Lambda(t) e^{z(t)}\right)
=e^{\frac{d^2}{2c_0}+\varep_0} \limsup_{t\to\infty} \left(\bar \Lambda(T_0+t) e^{\frac{d}{c_0}z_0(T_0+t)}\right)
=\ell e^{\frac{d^2}{2c_0}+\varep_0}.
\eeq 
It follows from \eqref{tJe} and \eqref{limy}  that 
\beq\label{ltJ}
\limsup_{k\to\infty}\widetilde J_k
\le \limsup_{k\in \N,k\to\infty}y(k+1)
\le \limsup_{t\to\infty}y(t)
\le \ell e^{\frac{d^2}{2c_0}+\varep_0}.
\eeq 
Combining this with \eqref{liJsum} gives 
\beq \label{lJk}
\limsup_{k\to\infty} J_k\le \ell e^{\frac{d^2}{2c_0}+\varep_0}.
\eeq 
Combining \eqref{wJk} with \eqref{lJk} yields 
\beqs
\limsup_{t\to\infty}\left( \max_{x\in\bar U} w^+(x,t)\right)\le \ell e^{\frac{d^2}{2c_0}+\varep_0}.
\eeqs
Passing $\varep_0\to 0$, we obtain \eqref{lel1}.

\medskip\noindent
\textbf{Step 8.} Now,  we consider, instead of \eqref{increase}, the case  when the function $q(t)$ in \eqref{increase} is decreasing.
By \eqref{newiden}, we have  $\Lambda(t)E(t+1)$ is decreasing.
With this property,  we estimate, instead of \eqref{uplamE}, 
\beqs
\Lambda(m)E_0(m+1)\le \int_{m-1}^{m} E_0(s+1) \Lambda(s) \d s.
\eeqs 
It implies
\begin{align*}
\widetilde J_k
&\le 
\frac{ \sum_{m=1}^k \int_{m-1}^{m} E_0(s+1) \Lambda(s) \d s }{E_0(k+1)}=
\frac{\int_{0}^{k} E_0(s+1) \Lambda(s) \d s }{E_0(k+1)}=Y(k),
\end{align*}
where
\beqs
Y(t)=\frac{\int_{0}^{t} E_0(s+1) \Lambda(s) \d s }{E_0(t+1)}.
\eeqs
We rewrite 
\beqs
Y(t)
=\frac{\int_{0}^{t} \exp\left(-\int_{s+1}^t e^{-z(x)} \d x\right) \Lambda(s) \d s}
{\exp\left(\int_{t}^{t+1} e^{-z(x)} \d x\right)}
=\frac{y(t)}{\exp\left(\int_{t}^{t+1} e^{-z(x)} \d x\right)}
\le y(t),
\eeqs
where $y(t)$ is defined by \eqref{ydef}.
Then, together with \eqref{limy},
\beqs
\limsup_{t\to\infty}Y(t)\le\limsup_{t\to\infty}y(t)\le \ell e^{\frac{d^2}{2c_0}+\varep_0}.
\eeqs
Thus,
\begin{align*}
\limsup_{k\to\infty}\widetilde J_k
\le \limsup_{k\in \N,k\to\infty}Y(k)
\le \limsup_{t\to\infty}Y(t)
\le \ell e^{\frac{d^2}{2c_0}+\varep_0},
\end{align*}
that is, we obtain \eqref{ltJ} again.
Using the same arguments as in Step 7 after \eqref{ltJ},  we obtain \eqref{lel1} again. The proof is complete.
\end{proof}

We obtain an asymptotic estimate for $|w(x,t)|$ next.

\begin{theorem}\label{NHthm2}
Suppose we have the same assumptions as in Theorem \ref{NHthm3} except that 
$|\widetilde Lw|$ replaces $\widetilde Lw$ in \eqref{Lwf} and $|w(x,t)|$ replaces $w^+(x,t)$ in \eqref{lamst}.
Then one has
\beq\label{lel2}
\limsup_{t\to\infty}\left( \max_{x\in\bar U} |w(x,t)|\right)\le e^\frac{d^2}{2c_0} \ell .
\eeq 
\end{theorem}
\begin{proof}
We apply Theorem \ref{NHthm3} to $w$ and $(-w)$, and then use inequality \eqref{absw} to obtain \eqref{lel2}.
We omit the details.
\end{proof}

Below, we have some examples for $z_0(t)$ and $\Lambda_*(t)$. Since $\Lambda_*(t)=\bar \Lambda(4t)$, it suffices to specify $\bar\Lambda$ instead of $\Lambda_*$.

\begin{proposition}\label{ex1}
Let $L>0$ and $\delta\ge 0$, and assume one of the following for sufficiently large $t>0$.

\begin{enumerate}[label=\tnum]

\item\label{case1} For any number $\gamma>0$, let $z_0(t)=\gamma\ln\ln t$ and $\bar \Lambda(t)=L (\ln t)^{-\gamma\left(\frac{d}{c_0}+\delta\right)}$. 

\item\label{case4}  
For any positive number $\gamma <c_0/d$ and number $\alpha \in\R$, let $z_0(t)=\gamma\ln ( t (\ln t)^\alpha)$  and 
$\bar \Lambda(t)=L [t(\ln t)^\alpha]^{- \gamma\left(\frac{d}{c_0}+\delta\right)}$.
In particular, if $\alpha=0$ then $z_0(t)=\gamma\ln t$  and $\bar \Lambda(t)=L t^{-\gamma\left(\frac{d}{c_0}+\delta\right)}$.

\item\label{case5} For any number $\alpha \in(-\infty,1)$, let $z_0(t)=\frac{c_0}{d}\ln ( t (\ln t)^\alpha)$   and 
$\bar \Lambda(t)=L [t(\ln t)^\alpha]^{-1-\frac{c_0\delta}{d}}$.
In particular, if $\alpha=0$ then  $z_0(t)=\frac{c_0}{d}\ln t$  and $\bar \Lambda(t)=L t^{-1-\frac{c_0\delta}{d}}$.
\end{enumerate}

Then the functions  $z_0(t)$ and $\bar \Lambda(t)$ satisfy \eqref{z2cond} and \eqref{monocond} for any $\varep>0$ with some sufficiently large numbers $T^*$ and $t_*$. Moreover,  the limit \eqref{Llim} holds with $\ell=L(1-\sign(\delta))$.
\end{proposition}
\begin{proof}
One can verify that condition \eqref{z2cond} is met in all the cases \ref{case1}, \ref{case4}, \ref{case5} for a certain sufficiently large $T^*$. We check the remaining condition \eqref{monocond} for some large $t_*$ next.
We construct the function $\bar\Lambda (t)$ in the following general form 
\beq \label{spebL}
\bar \Lambda(t)=L e^{-\left(\frac{d}{c_0}+\delta\right)z_0(t)}.
\eeq 
For $\varep>0$, let 
\beq \label{ELF}
E(t)=\exp\left(e^{-\frac{d^2}{2c_0}-\varep}\int_{t_*}^{t+1} e^{-\frac{d}{c_0}z_0(x)} \d x\right)
\text{ and } \mathcal{F}(t)=\bar \Lambda(t) E(t).
\eeq
It suffices to prove that $\mathcal{F}'(t)$ has one sign for sufficiently large $t$.
Taking the derivative of $\mathcal{F}(t)$ gives
\beq \label{genFp}
\begin{aligned}
    \mathcal{F}'(t)&=-\bar \Lambda(t) \left(\frac{d}{c_0}+\delta\right)z_0'(t) E(t)+\bar \Lambda(t) E(t)e^{-\frac{d^2}{2c_0}-\varep}e^{-\frac{d}{c_0}z_0(t+1)}\\
    &=\bar \Lambda(t)E(t)\left\{ -\left(\frac{d}{c_0}+\delta\right)z_0'(t) +e^{-\frac{d^2}{2c_0}-\varep}e^{-\frac{d}{c_0}z_0(t+1)}\right\}. 
\end{aligned}
\eeq 

\medskip
Case \ref{case1}. 
For sufficiently large $t_*$, the calculations  in \eqref{genFp} specifically yield 
\begin{align*}
\mathcal{F}'(t)
&=\bar \Lambda(t)E(t)\left\{ -\left(\frac{d}{c_0}+\delta\right) \frac{\gamma}{t\ln t} + \frac{e^{-\frac{d^2}{2c_0}-\varep}}{(\ln t)^{d\gamma/c_0}}  \right\}.
\end{align*}
From this, one can easily see that  $\mathcal{F}'(t)>0$ for large $t$.

\medskip
Case \ref{case4}. 
We have
$z_0'(t)=\frac\gamma t(1+\frac\alpha{\ln t}).$
Together with \eqref{genFp}, this implies
\begin{align*}
\mathcal{F}'(t)
&=\bar \Lambda(t)E(t)\left\{ -\frac{\gamma}{t}\left(\frac{d}{c_0}+\delta\right)\left(1+\frac\alpha{\ln t}\right)
+ \frac{e^{-\frac{d^2}{2c_0}-\varep}}{\left\{(t+1)(\ln(t+1))^\alpha\right\}^{d\gamma/c_0}}  \right\}.
\end{align*}
Since $d\gamma/c_0<1$, we $\mathcal{F}'(t)>0$ for large $t$.

\medskip
Case \ref{case5}. Taking $\gamma=c_0/d$, we have, same as in \ref{case4}, that
\beq \label{FFprime}
\mathcal{F}'(t)=\bar \Lambda(t)E(t)\left\{ -\frac{1}{t}\left(1+\gamma\delta\right)\left(1+\frac\alpha{\ln t}\right)
+ \frac{e^{-\frac{d^2}{2c_0}-\varep}}{(t+1)(\ln(t+1))^\alpha}   \right\}.
\eeq

If $\alpha\in(0,1)$, then the term $-(1+\gamma\delta)/t$ in \eqref{FFprime} dominates as $t\to\infty$, hence $\mathcal{F}'(t)<0$ for large $t$.

If $\alpha=0$, then \eqref{FFprime} becomes
\begin{align*}
\mathcal{F}'(t)
&=\bar \Lambda(t)E(t)\left\{ -\frac{1}{t}\left(1+\gamma\delta\right)
+ \frac{e^{-\frac{d^2}{2c_0}-\varep}}{t+1}  \right\}.
\end{align*}
Because $e^{-\frac{d^2}{2c_0}-\varep}<1\le 1+\gamma\delta$, it follows that $\mathcal{F}'(t)<0$ for large $t$.

If $\alpha<0$, then the last term in \eqref{FFprime} dominates, for large $t$, and hence $\mathcal{F}'(t)>0$. This completes the verification of condition \eqref{monocond} in Assumption \ref{AssumZ}.

\medskip
Now, with $\bar \Lambda(t)$ in \eqref{spebL}, it is obvious that  the limit superior  \eqref{Llim} is $\limsup_{t\to\infty}\left(L e^{-\delta z_0(t)}\right)$ which is $\ell=L$ when $\delta=0$, and is $\ell=0$ when $\delta>0$.
This proves the last statement of this proposition.
\end{proof}

Below, we give more examples for the case $\ell=0$ in \eqref{Llim}.
These aim to remove the dependence of the functions $\bar\Lambda(t)$ in Proposition \ref{ex1} on $\gamma$.

\begin{proposition}\label{ex2}
 For any numbers $\alpha,\beta,\gamma,L>0$, let  $z_0(t)=\gamma(\ln\ln t)^\alpha$ and $\bar \Lambda(t)=L t^{-\beta}$. Then the functions  $z_0(t)$ and $\bar \Lambda(t)$ satisfy \eqref{z2cond} and \eqref{monocond} for any $\varep>0$  with some sufficiently large numbers $T^*$ and $t_*$, and the limit \eqref{Llim} holds with $\ell=0$.
\end{proposition}
\begin{proof}
Let $E(t)$ and $\mathcal{F}(t)$ be as in \eqref{ELF}. Taking the derivative of $\mathcal{F}(t)$ gives
\begin{align*}
    \mathcal{F}'(t)&=-\beta \bar \Lambda(t) t^{-1} E(t)+\bar \Lambda(t) E(t)e^{-\frac{d^2}{2c_0}-\varep}e^{-\frac{d}{c_0}z_0(t+1)}
    =\bar \Lambda(t)E(t)\left\{ -\frac{\beta}t + \frac{e^{-\frac{d^2}{2c_0}-\varep}}{ e^{\frac{\gamma d}{c_0}[\ln \ln (t+1)]^\alpha}}\right\}. 
\end{align*}
For any $\varep'\in(0,1)$, one has, with sufficiently large $t$,  $e^{\frac{\gamma d}{c_0}[\ln \ln (t+1)]^\alpha}\le e^{\varep' \ln t}\le t^{\varep'}$.
Hence, $\mathcal{F}'(t)>0$ for large $t$. Also, by choosing $\varep'<\beta$, one induces
\beqs
\bar \Lambda(t) e^{\frac{d}{c_0}z_0(t)}=Lt^{-\beta} e^{\frac{\gamma d}{c_0}[\ln \ln t]^\alpha}\to 0\text{ as }t\to\infty.
\eeqs
This proves \eqref{Llim} with $\ell=0$.
\end{proof}

One can easily construct more complicated examples based on the above Propositions \ref{ex1} and \ref{ex2}.

\section{The nonlinear problem}\label{nonlinsec}
We now return to the nonlinear equation \eqref{maineq} and study the initial boundary value problem \eqref{nonIBVP}.
For a solution  $u\in C(\overline{U_\infty})\cap C_{x,t}^{2,1}(U_\infty)$ 
 of \eqref{nonIBVP}, Theorem \ref{maxprin2} and Corollary \ref{range} already provide some rough lower and upper bounds.
 To improve these estimates, we need to take advantage of the inhomogeneous  Growth Lemma, i.e., Lemma \ref{lemG2}. 
 It in turn requires a certain conversion of the nonlinear equation  \eqref{maineq} to a suitable linear problem.
 We will do just that in subsection \ref{CHsec} below.

\subsection{Reduction method}\label{CHsec}
Recall the nonlinear operator $L$ defined by \eqref{Ldef}. 
To remove the quadratic terms of the gradient in $Lu$, we use a transformation of the Bernstein--Cole--Hopf type \cite{Bernstein1912,Cole1951,Hopf1950}. 
In this subsection \ref{CHsec}, $T$ is an extended number in $(0,\infty]$, and $u$ is a function such that 
\beqs 
u\in C_{x,t}^{2,1}(U_T)\text{ and $u(U_T)\subset J$.}
\eeqs 

For any function $w\in C_{x,t}^{2,1}(U_T)$, we define  the function  $\mathcal L w$ on $U_T$  by 
\beq \label{reL}
\mathcal L w = w_t-\langle A(x,t),D^2w\rangle  +\widetilde B(x,t)\cdot \nabla w,\text{ where } \widetilde B(x,t)=B(x,t,u(x,t)). 
\eeq 
Clearly, $\mathcal L$ is a linear mapping from $C_{x,t}^{2,1}(U_T)$ to the space of functions defined on $U_T$.

The following lemma is an extension of \cite[Lemma 3.3]{HI3} to inhomogeneous problems.

\begin{lemma}\label{q-lin}

Below, for $\lambda\in \R$ and $C>0$, let $F_{\lambda,C}$ be any $C^2$-function on $J$ such that 
\beq \label{ftrans}
F'_{\lambda,C}(s)=C e^{\lambda P(s)} \text{ for any }  s\in J.
\eeq 
\begin{enumerate}[label=\tnum]    
    \item\label{Fsub} 
Assume there is a  constant $c_1\ge 0$ such that 
\beq\label{cond1}
 \xi^{\rm T} K(x,t)\xi \ge -c_1|\xi|^2 \text{ for all $(x,t)\in U_T$ and all $\xi\in\R^n$.} 
\eeq
For any numbers $\lambda\ge c_1/c_0$ and $C>0$, one has the function $w=F_{\lambda,C}(u)$ satisfies
\beq \label{Lw1}
\mathcal Lw\le F_{\lambda,C}'(u)Lu\text{  on }U_T.
\eeq 

 \item\label{Fsuper} 
Assume there is a constant $c_2\ge 0$ such that 
       \beq\label{cond2}
 \xi^{\rm T} K(x,t)\xi\le c_2|\xi|^2 \text{ for all $(x,t)\in U_T$ and all $\xi\in\R^n$.}
   \eeq 
For any numbers $\lambda\le  -c_2/c_0$ and $C>0$,  one has the function $w=F_{\lambda,C}(u)$  satisfies
\beq \label{Lw2}
\mathcal Lw\ge F_{\lambda,C}'(u)Lu\text{  on }U_T.
\eeq 
\end{enumerate}
\end{lemma}
\begin{proof}
Denoting $F=F_{\lambda,C}$, we have 
\beqs 
F'(s)>0 \text{ and } F''(s)=\lambda P'(s)F'(s)\text{ for any }s\in J.
\eeqs 
We obtain from \cite[Eq. (3.18)]{HI3} that 
\beq\label{LLrel}
\begin{aligned}
 \mathcal Lw
 &=   F'(u)\{ Lu- P'(u)(\nabla u)^{\rm T} K(x,t)\nabla u\} - F''(u) (\nabla u)^{\rm T} A(x,t)\nabla u\\
 &=   F'(u)Lu- P'(u)F'(u)(\nabla u)^{\rm T} K(x,t)\nabla u - F''(u) (\nabla u)^{\rm T} A(x,t)\nabla u.  
\end{aligned}
\eeq

\medskip
\ref{Fsub} In this case with $\lambda\ge c_1/c_0\ge 0$, one has  $F''(s)\ge 0$ and 
\beq\label{Fpos}
\begin{aligned}
  F''(s) \xi^{\rm T} A(x,t)\xi
  &\ge c_0F''(s)|\xi|^2=c_0\lambda P'(s)F'(s)|\xi|^2\\
  &\ge c_1 P'(s)F'(s)|\xi|^2
  \ge - P'(s)F'(s)\xi^{\rm T} K(x,t)\xi
\end{aligned}
\eeq 
for any $s\in J$ and $\xi\in\R^n$. The last inequality in \eqref{Fpos} comes from \eqref{cond1}.
Applying \eqref{Fpos} to $s=u(x,t)$, $\xi=\nabla u(x,t)$ and combining it with \eqref{LLrel},  we obtain
$ \mathcal Lw\le F'(u)Lu$
which proves \eqref{Lw1}.

\medskip
\ref{Fsuper} In this case $\lambda\le  -c_2/c_0\le 0$, one has $F''(s)\le 0$ and 
\beq\label{Fneg}
\begin{aligned}
  F''(s) \xi^{\rm T} A(x,t)\xi
  &\le c_0F''(s)|\xi|^2 = c_0\lambda P'(s)F'(s)|\xi|^2\\
  &\le -c_2 P'(s)F'(s)|\xi|^2
  \le - P'(s)F'(s)\xi^{\rm T} K(x,t)\xi
\end{aligned}
\eeq 
for any $s\in J$ and $\xi\in\R^n$.
Applying \eqref{Fneg} to $s=u(x,t)$, $\xi=\nabla u(x,t)$ and combining it with \eqref{LLrel}, we obtain
$ \mathcal Lw\ge F'(u)Lu$
which proves \eqref{Lw2}.
\end{proof}

Note that the function $F_{\lambda,C}$ in Lemma \ref{q-lin} exists and is strictly increasing on $J$.

\begin{example}\label{fcase}
We recall  the following examples  from \cite[Example 3.4]{HI3}.
\begin{enumerate}[label=\rnum]
    \item\label{Fe1} \emph{Case of isentropic, non-ideal gas flows.} With the choice in \eqref{isenchoice}, we can choose a particular function $F_{\lambda,C}$ as
   \beq\label{Fl1} F_\lambda(s)=\int_0^s e^{\lambda z^\gamma}dz \text{ for all } s\ge 0.
   \eeq 

    \item\label{Fe2}  \emph{Case of ideal gas.}
With the choice in \eqref{idealchoice}, we select $s_0=0$ in both cases of $J$.
For $\lambda=0$, we clearly can choose  $F_{\lambda,C}$ as
  \beqs
  F_\lambda(s)=s \text{ for all }s\in J.
  \eeqs 
For $\lambda\ne 0$, we can choose $F_{\lambda,C}$ as
  \beq\label{Fl21} 
  F_\lambda(s)=\frac1\lambda e^{\lambda s}\text{ for all }s\in J,\text{ or }
    F_\lambda(s)=\sign(\lambda)e^{\lambda s}\text{ for all }s\in J.
  \eeq 

    \item\label{Fe3}  \emph{Case of slightly compressible fluids.} 
We use the choice in \eqref{slighchoice}. When $\lambda\ne -1$, we can choose $F_{\lambda,C}$ as 
\beq\label{Fl3} 
F_\lambda(s)=\frac1{\lambda+1}s^{\lambda+1} \text{ for all } s>0,\text{ or }
F_\lambda(s)=\sign(\lambda+1)s^{\lambda+1} \text{ for all } s>0.
\eeq 
When $\lambda=-1$, we can choose $F_{\lambda,C}$ as  $F_\lambda(s)=\ln (s)$ for all $s>0$.
\end{enumerate}  
\end{example}

\subsection{Convergence of the solution}

We finally utilize all the techniques developed previously to study the original inhomogeneous, nonlinear problem \eqref{nonIBVP}. For different types of function $B(x,t,s)$ in \eqref{Ldef} and different scenarios of the forcing function and boundary data, we will establish the convergence of the solution, in the supremum norm  in $x$,  as time tends to infinity.

\begin{assumption}\label{lastA}
Let $A:U_\infty\to \mathcal M^{n\times n}_{{\rm sym}}$ and $K:U_\infty\to \mathcal M^{n\times n}$
be such that 
\begin{enumerate}[label=\tnum]
    \item $A(x,t)$ satisfies Assumption \ref{firstA} with $T=\infty$, and 
    \item  there are constants $c_1\ge 0$ and $c_2\ge 0$ such that 
    \beq\label{condall}
- c_1|\xi|^2\le  \xi^{\rm T} K(x,t)\xi\le c_2|\xi|^2 \text{ for all $(x,t)\in U_\infty$, all $\xi\in\R^n$.}
   \eeq 
\end{enumerate}
\end{assumption}

For the remainder of this section, let Assumption \ref{lastA} hold true. Let $u\in C(\overline{U_\infty})\cap C_{x,t}^{2,1}(U_\infty)$ be a solution of \eqref{nonIBVP}.
Let $\mathcal L$  be the linear mapping defined by \eqref{reL} with $T=\infty$.

Thanks to condition \eqref{condall} in Assumption \ref{lastA}, the matrix $K(x,t)$ satisfies  both \eqref{cond1} and \eqref{cond2} with  $T=\infty$.

\begin{assumption}\label{Pbound}
Either $J=\R$ or one of the following cases holds true.
\begin{enumerate}[label=\tnum]
\item\label{p1}  If $J=[a,\infty)$ or $J=(-\infty,b]$, then  $u(\overline{U_\infty})\subset J$.

\item\label{p2}  If $J=(a,\infty)$, then  $u(\overline{U_\infty})\subset [m,\infty)$ for some $m>a$.

\item\label{p3} If $J=(-\infty,b)$, then  $u(\overline{U_\infty})\subset (-\infty,M]$ for some $M<b$.

\item\label{p4} If $J$ is bounded, then  $u(\overline{U_\infty})\subset [m,M]$ with $m,M\in J$.
\end{enumerate}
\end{assumption}

\begin{remark}\label{ranrmk}
In the case of ideal gas flows \eqref{idealchoice} with $J=\R$, Assumption \ref{Pbound} trivially holds.
In general,  Assumption \ref{Pbound} is not meant to impose conditions on the solution $u(x,t)$ in $U_\infty$. Rather,  when needed, it requires some conditions only on the forcing term $f(x,t)$, the boundary data $g(x,t)$, and the initial data $u_0(x)$  by the virtue of  Corollary  \ref{range}.  
For example, concerning the case \ref{p2}, resp. \ref{p3},  we only need to require that the numbers $m_2$ and $\widetilde m_2$ in \eqref{mtm2}, resp., $m_1$ and $\widetilde m_1$ in \eqref{mtm1}, satisfy 
$m_2-\widetilde m_2>a$, resp., $m_1+\widetilde m_1<b$.    
\end{remark}

\begin{theorem}[$L_t^1$-forcing I]\label{NHL1}
In addition to Assumption \ref{Pbound}, assume the following.
\begin{enumerate}[label=\tnum]
\item\label{Hi} There are constants $c_B>0$ and $\gamma_0>0$ such that 
    \beq\label{B0}
    |B(x,t,s)|\le c_B(1+|s|)^{\gamma_0} \text{ for all $(x,t)\in U_\infty$ and $s\in J$.}
    \eeq  
\item\label{Hii} There is a continuous function $F:[0,\infty)\to[0,\infty)$ such that
    \beq\label{fc1}
    |f(x,t)|\le F(t) \text{ for all $ (x,t)\in U_\infty$.}
    \eeq 
\item\label{Hv} The function $F(t)$ in \ref{Hii} satisfies   
    \beq\label{fc2}
    \mu_0\eqdef \int_0^\infty F(\tau)\d \tau <\infty.
    \eeq 
\item\label{Hiii} There is $u_*\in J$ such that 
\beq\label{glim}
\lim_{t\to\infty}\left(\max_{x\in \Gamma}|g(x,t)-u_*|\right)=0.
\eeq
\end{enumerate}
Then one has
\beq\label{lim5}
\lim_{t\to\infty}\left(\max_{x\in \bar U}|u(x,t)-u_*|\right)=0.
\eeq
\end{theorem}
\begin{proof}
The proof is divided into five steps.

\medskip\noindent\textbf{Step 1.} 
By \eqref{glim}, $u$ is bounded on $\Gamma\times (0,\infty)$.
Thus, it is bounded on $\Gamma_\infty$ by some number $\mu_1>0$.
By Theorem \ref{maxprin2}\ref{MP3} and \eqref{fc2}, we have $u$ is bounded on $\overline{U_\infty}$ by
\beq\label{ub1}
|u(x,t)|\le \mu_2\eqdef \mu_1+\mu_0 \text{ for all } (x,t)\in \overline{U_\infty}.
\eeq 
By Assumption \ref{Pbound} and \eqref{ub1}, we have 
\beq \label{ranbar}
u(\overline{U_\infty})\subset [m_*,M_*]\text{ for some numbers $m_*,M_*\in J$.}
\eeq 
Let  $\mu_3$  be a positive number such that
\beq\label{Pmu}
|P(s)|\le \mu_3\text{ for all }s\in[m_*,M_*].
\eeq
Then 
\beq\label{Pub1}
|P(u(x,t))|\le \mu_3 \text{ for all } (x,t)\in \overline{U_\infty}.
\eeq
Thanks to properties \eqref{B0} and \eqref{ub1}, we have $\widetilde B(x,t)$ is bounded in $U_\infty$ by
\beq\label{tBbo1}
|\widetilde B(x,t)|\le \mu_4\eqdef c_B(1+\mu_2)^{\gamma_0} \text{ for all } (x,t)\in U_\infty.
\eeq

\medskip\noindent\textbf{Step 2.} 
Let 
\beqs 
\lambda_1>0\text{ and } \lambda_2<0\text{ such that }
\lambda_1\ge c_1/c_0\text{ and } \lambda_2\le -c_2/c_0.
\eeqs
For $j=1,2$, let $F_{\lambda_j}$ be the function $F_{\lambda_j,1}$ defined by \eqref{ftrans}, and define, thanks to \eqref{ranbar},  
\beqs 
w_j=F_{\lambda_j}(u)\text{ on }\overline{U_\infty}.
\eeqs 
 By Lemma \ref{q-lin} applied to $T=\infty$, $w:=w_1$ in part \ref{Fsub} and $w:=w_2$ in part \ref{Fsuper}, we have 
 \beq \label{Lsign}
 \mathcal Lw_1\le e^{\lambda_1 P(u(x,t))}f(x,t)\text{ and }\mathcal Lw_2\ge e^{\lambda_2 P(u(x,t))}f(x,t)\text{  in }U_\infty.
 \eeq 
Denote $C_1=e^{\mu_3\max\{\lambda_1,-\lambda_2\}}$. Then combining \eqref{Pub1} with \eqref{fc1} gives
\beq \label{elamP}
 e^{\lambda_1 P(u(x,t))}|f(x,t)|\le \mathcal{F}(t):= C_1 F(t)
 \text{ and }
 e^{\lambda_2 P(u(x,t))}|f(x,t)|\le \mathcal{F}(t)\text{  in }U_\infty.
 \eeq
 For $j=1,2$, define $w_{*,j}=F_{\lambda_j}(u_*)$ and 
$\bar w_j=w_j-w_{*,j}$ on $\overline{U_\infty}$.
From \eqref{Lsign} and \eqref{elamP} one has
 \beq\label{bLw1}
\mathcal L\bar w_1= \mathcal Lw_1\le \mathcal{F}(t)
\text{ and }
\mathcal L\bar w_2= \mathcal Lw_2\ge -\mathcal{F}(t) \text{  in }U_\infty.
 \eeq 

 \medskip\noindent\textbf{Step 3.} 
For $j=1,2$, by the local Lipschitz continuity of $F_{\lambda_j}(s)$ for $s\in J$ near $u_*$ and thanks to \eqref{glim}, there is a positive constant $C_*$ such that, for sufficiently large $t$ and all $x\in \Gamma$, one has
\beq\label{wFg}
|w_j(x,t)-w_{*,j}|=|F_{\lambda_j}(g(x,t))-F_{\lambda_j}(u_*)|
\le C_* |g(x,t)-u_*|.
\eeq
Together with \eqref{glim} again, this implies, for $j=1,2$,
\beq \label{wblim}
\lim_{t\to\infty}\left(\max_{x\in\Gamma} |\bar w_j(x,t)|\right)=0.
\eeq

 Using the first inequality in \eqref{bLw1}, we can apply Theorem \ref{NHthm}\ref{dec1} to operator $\widetilde L:=\mathcal L$, function $w:=\bar w_1$, number $M_2:=\mu_4$ and function $F_1(t):=\mathcal{F}(t)$. 
 We obtain from \eqref{lsup1} that
\beq\label{liw5}
\limsup_{t\to\infty} \left[\max_{x\in\bar U}\bar  w_1^+(x,t)\right]
\le \frac{2-\eta_*}{1-\eta_*}\left(
\limsup_{t\to\infty}\left[\max_{x\in\Gamma} \bar w_1^+(x,t)\right]
+\limsup_{t\to\infty}\int_t^{t+T_*} \mathcal{F}(\tau)\d\tau 
\right).
\eeq
For the last limit superior, it follows from condition \eqref{fc2} that 
\beq\label{lifcond}
\lim_{t\to\infty} \int_t^{t+T_*}\mathcal{F}(\tau)\d \tau=C_1 \lim_{t\to\infty} \int_t^{t+T_*}F(\tau)\d \tau =0.
\eeq
Then, by combining \eqref{liw5} with the limits \eqref{wblim} and \eqref{lifcond}, one obtains
\beq\label{lim1}
\limsup_{t\to\infty} \left[\max_{x\in\bar U} \bar w_1^+(x,t)\right]
= 0.
\eeq

Similarly,  using the second inequality in \eqref{bLw1}, we can apply Theorem \ref{NHthm}\ref{dec1} to operator  $\widetilde L=\mathcal L$, function $w:=-\bar w_2$, number $M_2=\mu_4$ and function $F_1(t):=\mathcal{F}(t)$. We obtain from \eqref{lsup1} that
\beqs
\limsup_{t\to\infty} \left[\max_{x\in\bar U}\bar w_2^-(x,t)\right]
\le \frac{2-\eta_*}{1-\eta_*}\left(
\limsup_{t\to\infty}\left[\max_{x\in\Gamma}  \bar w_2^-(x,t)\right]
+\limsup_{t\to\infty}\int_t^{t+T_*} F(\tau)\d\tau 
\right).
\eeqs
Thus, thanks to \eqref{wblim} and \eqref{lifcond} again,
\beq\label{lim2}
\limsup_{t\to\infty} \left[\max_{x\in\bar U}\bar w_2^-(x,t)\right]
= 0.
\eeq

\medskip\noindent\textbf{Step 4.}    
Define $\bar u=u-u_*$ on $\overline{U_\infty}$.
Denote $C_2=e^{\mu_3\min\{-\lambda_1,\lambda_2\}}>0$.
For $j=1,2$, we have from \eqref{ftrans} and \eqref{Pmu} that
$$F_{\lambda_j}'(s)=e^{\lambda_j P(s)}\ge C_2 \text{ for }  s\in[m_*,M_*].$$
By the Mean Value Theorem, one has,  for $j=1,2$,
\begin{align} \label{suF}
C_2(s-u_*)&\le F_{\lambda_j}(s)-F_{\lambda_j}(u_*) \text{ for } s\in[u_*,M_*],\\
\label{Fsu}
C_2(u_*-s)&\le F_{\lambda_j}(u_*)-F_{\lambda_j}(s) \text{ for } s\in[m_*,u_*).
\end{align}
When  $j=1$, we have, based on  \eqref{suF}, that
\beq\label{uF1} 
(s-u_*)^+\le C_2^{-1}(F_{\lambda_1}(s)-F_{\lambda_1}(u_*))^+
\text{ for all }s\in[m_*,M_*].
\eeq
(For, \eqref{uF2} clearly holds true when $m_*\le s\le u_*$.)
Similarly, when  $j=2$, we have, based on  \eqref{Fsu}, that
\beq\label{uF2}
(s-u_*)^-
\le C_2^{-1}(F_{\lambda_2}(s)-F_{\lambda_2}(u_*))^-\text{ for all } s\in[m_*,M_*].
\eeq 
Applying \eqref{uF1} and \eqref{uF2} to $s=u(x,t)$ yields
\beq\label{upmw}
\bar u^+(x,t)\le C_2^{-1}w_1^+(x,t)\text{ and }\bar u^-(x,t)\le C_2^{-1}\bar w_2^-(x,t).
\eeq 

\medskip\noindent\textbf{Step 5.} 
Thanks to inequalities \eqref{absw}, \eqref{upmw} and the limits \eqref{lim1}, \eqref{lim2}, we obtain 
\begin{align*}
&\limsup_{t\to\infty}\left[\max_{x\in \bar U}|\bar u(x,t)|\right]\le \max\left\{ \limsup_{t\to\infty}\left[\max_{x\in \bar U}\bar u^+(x,t)\right] ,
\limsup_{t\to\infty}\left[\max_{x\in \bar U}\bar u^-(x,t)\right] \right\}\\
&\le  C_2^{-1}\max\left \{\limsup_{t\to\infty}\left[\max_{x\in \bar U} \bar w_1^+(x,t)\right],
\limsup_{t\to\infty}\left[\max_{x\in \bar U}\bar w_2^-(x,t)\right]\right\}= 0.
\end{align*}
Thus, we obtain \eqref{lim5}.
\end{proof}

Next, we consider the case when $B(x,t,s)$ may not be uniformly bounded in $t\in[0,\infty)$.
More specifically, the constant $c_B$ in \eqref{B0} is replaced with an unbounded  time-dependent function. For the sake of having concise statements below, we say that two nonnegative functions $z_0(t)$ and $\bar \Lambda(t)$, which are defined for sufficiently large $t$, \emph{satisfy Assumption \ref{AssumZ}} if the conditions from \eqref{z2cond} to \eqref{monocond} are met for some numbers $T^*$, $t_*$ and $\varep_*$.

\begin{theorem}[$L_t^1$-forcing II]\label{NHL2}
In addition to  Assumption \ref{Pbound}, we assume the following.
\begin{enumerate}[label=\tnum]
    \item\label{Ni} There is as continuous, increasing function $\bar b:[0,\infty)\to[0.\infty)$ which tends to infinity as $t\to\infty$, and a number $\gamma_0>0$ such that
\beq\label{B3}
|B(x,t,s)|\le \bar b(t)(1+|s|)^{\gamma_0}\text{ for all } (x,t)\in U_\infty\text{ and } s\in J.
\eeq 
    \item\label{Nii} There is a continuous function $F:[0,\infty)\to[0,\infty)$ that satisfies \eqref{fc1} and \eqref{fc2}.
    \item\label{Niii} There are numbers $u_*\in J$, $T^*>0$ and a decreasing, continuous function $\widetilde\Lambda(t)\ge 0$ on $[T^*,\infty)$ such that
\beq\label{lamnon}
\max_{(x,\tau)\in \Gamma\times[t,t+1]}\left|g(x,\tau)-u_* \right|
+\int_t^{t+1}F(\tau)\d \tau 
\le \widetilde \Lambda(t) \text{ for all }t\ge T^*.
\eeq
    \item\label{Niv} For any $\mu>0$, the functions  $z_0(t):=\mu \bar b (t)$ and  $\bar\Lambda(t):= \widetilde \Lambda(t/4)$  satisfy Assumption \ref{AssumZ} and the limit \eqref{Llim} holds with $\ell:=\ell_\mu=0$.
\end{enumerate}
Then one has \eqref{lim5}.
\end{theorem}
\begin{proof}
We follow the proof of Theorem \ref{NHL1} step by step.

\noindent\textbf{Step 1.} 
 We proceed up to  \eqref{Pub1}, and then replace \eqref{tBbo1} with the following, thanks to \eqref{B3} and \eqref{ub1},
\beq\label{tBbo2}
|\widetilde B(x,t)|\le \widehat z_0(t):=\mu_6 \bar b(t)\text{ with }\mu_6=(1+\mu_2)^{\gamma_0}.
\eeq 

\noindent\textbf{Step 2.} Unchanged.

\noindent\textbf{Step 3.}  We proceed up to \eqref{wblim}.
Let $\mu_7=\max\{C_*,C_1\}$. 
For $j=1,2$, using \eqref{wFg}, the definition of $\mathcal{F}(t)$ in \eqref{elamP}, 
and the assumption \eqref{lamnon} with function  $\widetilde \Lambda(t)$, we can select a sufficiently large $T_*>0$ and have 
\beq\label{newgFL}
\max_{\Gamma\times[t,t+1]}|\bar w_j(x,t)|
+\int_t^{t+1}\mathcal{F}(\tau)\d \tau
\le \widehat \Lambda(t):= \mu_7 \widetilde \Lambda(t) \text{ for all }t\ge T^*.
\eeq

For any numbers $\mu>0$ and $\gamma>0$, by assumption \ref{Niv},  the functions  $z_0(t):=\mu \bar b(t)$ and $\bar \Lambda(t):=\gamma \widetilde \Lambda(t/4)$ satisfy Assumption \ref{AssumZ} and the limit \eqref{Llim} holds with 
\beqs
\ell:=\ell_{\mu,\gamma}=\limsup_{t\to\infty} \left(\gamma\widetilde \Lambda(t/4) e^{\frac{d}{c_0}\mu z_0(t)}\right)
=\gamma\limsup_{t\to\infty} \left(\widetilde \Lambda(t/4) e^{\frac{d}{c_0}\mu \bar b(t)}\right)=\gamma\ell_\mu
=0.
\eeqs 
The last identity is, of course, due to the last statement in \ref{Niv}.

(a) We consider the first inequality in \eqref{bLw1} with $\bar w_1$. 
By the remark right after \eqref{newgFL},   the functions  $z_0(t):=\widehat z_0(t)=\mu_6 \bar b(t)$ and $\bar \Lambda(t):=\widehat \Lambda(t/4)=\mu_7 \widetilde \Lambda(t/4)$ satisfy Assumption \ref{AssumZ} and the limit \eqref{Llim} holds with $\ell:=\ell_{\mu_6,\mu_7}=0$. 
We then apply Theorem \ref{NHthm3} to $w:=\bar w_1$, $F(t):=\mathcal{F}(t)$, $z_0(t):=\widehat z_0(t)$ in \eqref{tBbo2}, $\bar \Lambda(t):=\widehat \Lambda(t/4)$ in \eqref{newgFL} and $\ell:=\ell_{\mu_6,\mu_7}=0$.
Thus, we have from \eqref{lel1} that 
\beqs 
\limsup_{t\to\infty} \left[\max_{x\in \bar U} \bar w_1^+(x,t)\right]\le  e^\frac{d^2}{2c_0}\ell_{\mu_6,\mu_7}=0,
\eeqs 
that is, we obtain \eqref{lim1} again.

(b) Consider the second inequality in \eqref{bLw1} with $\bar w_2$.
Same as part (a), by applying Theorem \ref{NHthm3}  to $w:=-\bar w_2$, $F(t):=\mathcal{F}(t)$, $z_0(t):=\widehat z_0(t)$  and $\bar \Lambda(t):=\widehat \Lambda(t/4)$, we have 
\beqs  
\limsup_{t\to\infty}\left[ \max_{x\in\bar U}\bar w_2^-(x,t)\right]\le e^\frac{d^2}{2c_0}\ell_{\mu_6,\mu_7}=0.
\eeqs 
Thus, we obtain \eqref{lim2} again.

\noindent\textbf{Steps 4 and 5.} 
Unchanged. At the end, we obtain the limit \eqref{lim5}.
\end{proof}

\begin{example}\label{ex3}
Assume, for sufficiently large $t$,
    \beqs
    \bar b(t)=N_1 (\ln\ln t)^\alpha,\quad F(t)=N_2 t^{-\beta_0},\quad \max_{x\in\Gamma}\left|g(x,t)-u_*\right|\le N_3 t^{-\beta_1} ,
    \eeqs
    where $\alpha>0$, $\beta_0>1$, $\beta_1>0$ and $N_1,N_2,N_3>0$. Then the left-hand side of \eqref{lamnon} is bounded by 
    \beqs
    N_3 t^{-\beta_1}+\frac{N_2}{\beta_0-1}t^{1-\beta_0}\le \widetilde\Lambda(t):=N_4 t^{-\beta_2},
    \eeqs
    where $\beta_2=\min\{\beta_1,\beta_0-1\}$ and $N_4>0$. 
By the virtue of  Proposition \ref{ex2}, the condition \ref{Niv} in Theorem \ref{NHL2} is met. 
\end{example}

\begin{theorem}[Non-$L^1_t$ forcing I]\label{NHL3}
Consider  
\beq\label{ssca} 
J=(0,\infty)\text{ and $P(s)=\ln s$ for all $s\in J$.}
\eeq 
Assume \ref{Ni} in Theorem \ref{NHL2} and \ref{Hii} in Theorem \ref{NHL1} and 
   \beq\label{fc3}
 \int_0^\infty F(\tau)\d \tau =\infty.
    \eeq 
Suppose there exist  a number $u_*\in J$, a positive number $\lambda_1\ge c_1/c_0$, 
and two continuous  functions $\mathcal F(t)$ and  $\widetilde \Lambda(t)$ defined for sufficiently large $t$  such that the following are satisfied.
\begin{enumerate}[label=\tnum] 
\item\label{Lzo} The function $\mathcal F(t)$ is increasing and, for sufficiently large $t$,
\beq\label{iFmF}
\int_0^t F(\tau)\d \tau\le \mathcal{F}(t).
\eeq

\item\label{Li} The function $\widetilde \Lambda(t)$ is decreasing and, for sufficiently large $t$,  
    \beq\label{newLam}
\max_{(x,\tau)\in \Gamma\times [t,t+1]} |g(x,\tau)-u_*|
+\int_t^{t+1}F(\tau) \mathcal{F}(\tau)^{\lambda_1}\d\tau
    \le  \widetilde \Lambda(t).
    \eeq
    
\item\label{Liii} For any $\mu>0$, the functions $z_0(t):= \mu \bar b(t) \mathcal{F}(t)^{\gamma_0}$ and 
  $\bar\Lambda(t):=  \widetilde \Lambda(t/4)$ satisfy Assumption \ref{AssumZ} and the limit \eqref{Llim} holds with $\ell:=\ell_\mu =0$.
\end{enumerate} 

If, additionally,  
\beq \label{Liv} 
u(\overline{U_\infty})\subset [m,\infty)\text{ for some number $m>0$,}
\eeq 
then one has \eqref{lim5}.
\end{theorem}
\begin{proof}
We denote 
    \beq\label{Lg} 
    \Lambda_g(t)= \max_{(x,\tau)\in\Gamma\times [t,t+1]} |g(x,\tau)-u_*|.
    \eeq 

\medskip\noindent 
\textbf{Step 1.} 
Thanks to \eqref{fc3} and \eqref{iFmF}, one has 
\beq\label{unb} 
\mathcal{F}(t)\to\infty\text{ as }t\to\infty.
\eeq 
According to Condition \ref{Liii} with $\mu=1$,  we have the limit 
\beq\label{copylim}
\lim_{t\to\infty} \left\{\widetilde \Lambda(t/4) \exp\left (\frac{d}{c_0}\bar b(t) \mathcal{F}(t)^{\gamma_0}\right)\right\}=0.
\eeq
By the fact $\lim_{t\to\infty} \bar b(t)=\infty$ 
in Condition \ref{Ni} of Theorem \ref{NHL2} and  
and property \eqref{unb} of $\mathcal F(t)$, we have from the limit  \eqref{copylim} that $\widetilde \Lambda(t)\to 0$ as $t\to\infty$.
This and \eqref{newLam} in turn imply
\beq\label{gclous}
\Lambda_g(t)\to 0\text{ as }t\to\infty.
    \eeq
As a consequence of \eqref{gclous}, one has 
\beq \label{uGam}
\text{$u$ is bounded on $\Gamma_\infty$.}
\eeq 
By estimate \eqref{maxabs} in Theorem \ref{maxprin2}\ref{MP3}, taking into account \eqref{uGam}  and properties  \eqref{iFmF}, \eqref{unb}, we can estimate, for $(x,t)\in \overline{U_\infty}$ with large $t$, 
\beq\label{umuF}
|u(x,t)|\le \mu_1 \mathcal{F}(t), \text{ where $\mu_1$ is some positive number.}
\eeq
Together with \eqref{B0}, this implies, for $(x,t)\in U_\infty$ with large $t$,
\beq\label{tBF1}
|\widetilde B(x,t)|\le \widehat z_0(t):=\mu_2 \bar b(t) \mathcal{F}(t)^{\gamma_0}
\text{ for some number $\mu_2>0$.}
\eeq
 
\medskip\noindent 
\textbf{Step 2.} 
Let $\lambda_2<0$ such that $-1\ne \lambda_2\le -c_2/c_0$.
For $j=1,2$, define the function  $w_j=F_{\lambda_j}(u):=\frac{1}{\lambda_j+1}u^{\lambda_j+1}$ on $\overline{U_\infty}$, see the first formula in \eqref{Fl3} of Example \ref{fcase}\ref{Fe3}. 
Let $\bar u=u-u_*$, and, for $j=1,2$, let $w_{*,j}=F_{\lambda_j}(u_*)$ and define $\bar w_j=w_j-w_{*,j}$ on $\overline{U_\infty}$.
Since property \eqref{gclous} implies \eqref{glim},  we obtain inequality  \eqref{wFg} again, i.e., for all $x\in\Gamma$ and sufficiently large $t$,
    \beqs 
|\bar w_j(x,t)|\le C_* |g(x,t)-u_*| \text{ for $j=1,2$.}
\eeqs
Thus, for large $t$ and $j=1,2$,
\beq\label{wbc3}
\max_{\Gamma\times [t,t+1]}|\bar w_j|\le C_* \max_{\Gamma\times [t,t+1]}|g-u_*|
= C_* \Lambda_g(t).
\eeq

\medskip\noindent 
\textbf{Step 3.} 
According to Lemma \ref{q-lin}\ref{Fsub} with $\lambda=\lambda_1$ and $C=1$, one has, for $(x,t)\in U_\infty$, that 
\beq \label{LwFp}
\mathcal L w_1(x,t)\le f(x,t) F_{\lambda_1}'(u(x,t))=f(x,t) u(x,t)^{\lambda_1}.
\eeq 
We  use \eqref{fc1} to estimate $|f(x,t)|$ and use \eqref{umuF} to estimate $|u(x,t)|$ in \eqref{LwFp}.
Then, for $(x,t)\in U_\infty$,
\beq \label{cF1def}
\mathcal L \bar w_1 =\mathcal L w_1 
\le \mathcal{F}_1(t):=\mu_3 F(t) \mathcal{F}^{\lambda_1}(t), \text{ where }\mu_3=\mu_1^{\lambda_1}.
\eeq 
By \eqref{wbc3}, the definition of $\mathcal{F}_1(t)$ in \eqref{cF1def}, and inequality \eqref{newLam} we have, for large $t$,  
\beqs
\max_{\Gamma\times[t,t+1]}|\bar w_1|
+\int_t^{t+1}\mathcal{F}_1(\tau)\d\tau
\le 
C_* \Lambda_g(t)
+\mu_3 \int_t^{t+1}F(\tau) \mathcal{F}(\tau)^{\lambda_1}\d\tau 
\le \mu_4 \widetilde \Lambda(t),
\eeqs 
where $\mu_4=\max\{C_*,\mu_3\}$.
Because of Condition \ref{Liii},  by using the same argument as in Step 3 of the proof of Theorem \ref{NHL2} above, we can apply Theorem \ref{NHthm3} to $w:=\bar w_1$, $F(t):=\mathcal{F}_1(t)$, $z_0(t):=\widehat z_0(t)$ in \eqref{tBF1}, $\bar \Lambda(t):=\mu_4 \widetilde \Lambda(t/4)$ and 
\beq\label{Lmmlim}
\ell:=\ell_{\mu_2,\mu_4}
\eqdef \limsup_{t\to\infty} \left(\mu_4\widetilde \Lambda(t/4) e^{\frac{d}{c_0}\widehat z_0(t)}\right)
=\mu_4\limsup_{t\to\infty} \left(\widetilde \Lambda(t/4) e^{\frac{d}{c_0}\mu_2\bar b(t) \mathcal{F}(t)^{\gamma_0}}\right)=\mu_4\ell_{\mu_2}
=0.
\eeq 
It follows from inequality \eqref{lel1} that 
\beq \label{liw6}
\limsup_{t\to\infty} \left[\max_{x\in \bar U}\bar w_1^+(x,t)\right]\le \mu_4 e^\frac{d^2}{2c_0}\ell=0.
\eeq
Since $\bar w_1\le \bar w_1^+$, we obtain
\beqs
\limsup_{t\to\infty} \left[\max_{x\in \bar U}\bar w_1(x,t)\right]\le 0.
\eeqs
Consequently,
\beqs
\limsup_{t\to\infty}\left[\max_{\bar U} w_1(x,t) \right] \le w_{1,*},
\text{ that is, }
\limsup_{t\to\infty}\left[\max_{\bar U} F_{\lambda_1}(u(x,t)) \right]
\le F_{\lambda_1}(u_*).
\eeqs
From this we find that, see e.g. \cite[Inequalities  (4.60)  and (4.61)]{HI3},
\beq\label{ue2}
\limsup_{t\to\infty}\left[\max_{x\in \bar U} u(x,t) \right] \le u_*.
\eeq
Consequently,
\beq\label{ue3}
\limsup_{t\to\infty}\left[\max_{x\in \bar U} \bar u(x,t) \right] \le 0.
\eeq

\medskip\noindent 
\textbf{Step 4.}  
By \eqref{Liv}, we have $u(x,t)\ge m>0$ for all $(x,t)\in \overline{U_\infty}$.
According to Lemma \ref{q-lin}\ref{Fsuper} with $\lambda=\lambda_2<0$ and $C=1$, one has, for $(x,t)\in U_\infty$, that 
\beqs 
\mathcal L w_2\ge f(x,t) u^{\lambda_2}\ge -F(t)u^{\lambda_2}
\ge -F(t)m^{\lambda_2}=-\mathcal{F}_2(t),\text{ where } \mathcal{F}_2(t)=F(t)m^{\lambda_2}.
\eeqs 
This implies 
\beqs 
\mathcal L \bar w_2 =\mathcal L w_2 \ge -\mathcal{F}_2(t).
\eeqs 
Consider sufficiently large $t$.
By \eqref{unb}, we have $m^{\lambda_2}\le \mathcal{F}(t)^{\lambda_1}$.
Together with \eqref{wbc3} and \eqref{newLam}, this yields
\beqs
\max_{\Gamma\times[t,t+1]}|\bar w_2|
+\int_t^{t+1}\mathcal{F}_2(\tau)\d\tau
\le C_* \Lambda_g(t)
+\int_t^{t+1}F(\tau) \mathcal{F}(\tau)^{\lambda_1}\d\tau 
\le \mu_5 \widetilde \Lambda(t),    
\eeqs
 where $\mu_5=\max\{C_*,1\}$.
Applying Theorem \ref{NHthm3} to $w:=-\bar w_2$,  $F(t):=\mathcal{F}_2(t)$, $z_0(t):=\widehat z_0(t)$ in \eqref{tBF1}, $\bar \Lambda(t):=\mu_5 \widetilde \Lambda(t/4)$ and, same as \eqref{Lmmlim},  $\ell:=\ell_{\mu_2,\mu_5}=\mu_5\ell_{\mu_2}=0$, we obtain from 
\eqref{lel1} that
\beqs 
\limsup_{t\to\infty}\left[\max_{x\in\bar U}(-\bar w_2(x,t))^+\right]\le e^\frac{d^2}{2c_0}\ell_{\mu_2,\mu_5}=0.
\eeqs
Hence,
\beqs  
0\ge \limsup_{t\to\infty}\left[\max_{x\in\bar U}(-\bar w_2(x,t))\right]
= -\liminf_{t\to\infty}\left[\min_{x\in\bar U}\bar w_2(x,t)\right].
\eeqs 
This implies
\beqs  
\liminf_{t\to\infty}\left[\min_{x\in\bar U}\bar w_2(x,t)\right]\ge 0.
\eeqs 
Therefore, one has, similar to \eqref{ue2}, 
\beqs
\liminf_{t\to\infty} \left[\min_{x\in \bar U} u(x,t) \right]\ge  u_*,
\eeqs 
and, consequently,
\beq\label{ue5}
\liminf_{t\to\infty} \left[\min_{x\in \bar U} \bar u(x,t) \right]\ge  0.
\eeq

\medskip\noindent 
\textbf{Step 5.} Let $\varep>0$ be arbitrary.
By \eqref{ue3} and \eqref{ue5}, we have for sufficiently large $t$ that
\beqs
-\varep \le \min_{x\in \bar U} \bar u(x,t)\le \bar u(x,t) \le \max_{x\in \bar U} \bar u(x,t) \le \varep \text{ for all }x\in \bar U,
\eeqs
which implies $|\bar u(x,t)|\le \varep$  for all $x\in \bar U$.
This proves the limit \eqref{lim5}.
\end{proof}

Same as Remark \ref{ranrmk}, we comment again that condition \eqref{Liv} in Theorem \ref{NHL3} can be met thanks purely to the initial and boundary data. This comes from Corollary \ref{range} with the condition $m_2>\widetilde m_2$. In that case, $m=m_2-\widetilde m_2>0$.

\begin{example}\label{ex4}
We give an example for the functions $\bar b(t)$, $F(t)$ and $\Lambda_g(t)$, see \eqref{Lg}, so that Theorem \ref{NHL3} can apply.
    Let $\delta$ and $\beta$ be any numbers in the interval $(0,1)$. 
    Assume, for sufficiently large $t$,
\beqs 
\bar b(t)=N_1 (\ln\ln t)^{\gamma_1},\quad F(t)= \frac{N_2}{t (\ln t) (\ln\ln t)^\delta},
\quad \Lambda_g(t)\le  \frac{N_3}{t^\beta}
\eeqs 
for some positive numbers $\gamma_1$, $N_1$, $N_2$, $N_3$. Below, $N_j$, for $4\le j\le 7$, are appropriate positive numbers and $t$ is sufficiently large.
We estimate
\beqs
 \int_0^t F(\tau)\d \tau \le \mathcal{F}(t):=N_4(\ln\ln t)^{1-\delta}.
\eeqs
Because $F(t)$ is decreasing and $\mathcal{F}(t)$ is increasing, one has
\beqs
 \Lambda_g(t)
+\int_t^{t+1}F(\tau) \mathcal{F}(\tau)^{\lambda_1}\d\tau 
\le 
\Lambda_g(t)
+F(t) \mathcal{F}(t+1)^{\lambda_1}
\le  \widetilde \Lambda(t):=N_5 t^{-\beta}.
\eeqs
Let $\mu>0$ be arbitrary. Set
\beqs 
z_0(t)=\mu \bar b(t) \mathcal{F}(t)^{\gamma_0}=N_6 (\ln\ln t)^{\gamma_2},\text{ where }
\gamma_2=\gamma_1+\gamma_0(1-\delta)>0.
\eeqs 
By Proposition \ref{ex2}, the functions $z_0(t)$ and  $\bar\Lambda(t):= \widetilde \Lambda(t/4)= N_7 t^{-\beta}$ satisfy Assumption \ref{AssumZ} and the limit \eqref{Llim} holds with $\ell=0$.
\end{example}

The next result shows that the convergence in Theorem \ref{NHL3} can still hold true even when  $u_*\not \in J$.

\begin{theorem}[Non-$L^1_t$ forcing II]\label{NHL4}
Assume the same as in Theorem \ref{NHL3} with $u_*=0$, except for the last requirement \eqref{Liv}. 
Then one has
\beq\label{mainlim2}
\lim_{t\to\infty}\left\{ \max_{x\in \bar U}u(x,t)\right\} =0.
\eeq
\end{theorem}
\begin{proof}
    We follow the proof of Theorem \ref{NHL3}. In this case $g\ge 0$, $u\ge 0$ and $u_*=0$.
    Although $P(s)$ is not defined at $s=u_*=0$, we have for $\lambda_1>0$ the transformation $F_{\lambda_1}(s)=\frac{s^{\lambda_1+1}}{\lambda_1+1}$ for $s>0$ can be extended to the function
    $\widehat F(s)=\frac{s^{\lambda_1+1}}{\lambda_1+1}$  for all $s\ge 0$.
    Note that $\widehat F'(s)=s^{\lambda_1}$ is bounded near $0$. In the proof of Theorem \ref{NHL3}, we replace $F_{\lambda_1}$ with $\widehat F$ and proceed up to \eqref{liw6}. In this case $\bar w_1=w_1=\widehat F(u)\ge 0$. Then inequality \eqref{liw6}  reads as
\beqs
\limsup_{t\to\infty} \left[\frac1{\lambda_1+1}\max_{x\in \bar U} u^{\lambda_1+1}(x,t)\right]\le 0.
\eeqs
This implies, thanks to the fact $u\ge 0$, the desired limit \eqref{mainlim2}.
\end{proof}

\begin{remark}
The following final remarks are in order.
\begin{enumerate}[label=\rnum]
\item  Same as Theorem \ref{NHL4}, one can consider \eqref{ssca}  and $u_*=0$ for the $L_t^1$-forcing cases such as in Theorems \ref{NHL1} and \ref{NHL2}.

\item  For more cases of $u(\overline{U_\infty})\subset \bar J$ and/or $u_*\in \partial J\setminus J$, see  \cite[Theorem 4.7]{HI3}.

\item It is open that whether Theorem \ref{NHL3} can be established for isentropic fluid flows, see \eqref{isenchoice} and \eqref{idealchoice}. The main obstacle is that the term $F_{\lambda_1}'(u(x,t))$ in \eqref{LwFp} can grow too large as $t\to\infty$, see \eqref{Fl1} and \eqref{Fl21}, and, hence, Theorem \ref{NHthm3} may not apply. 

\item In this paper, we studied only the Darcy law \eqref{Darcy} for single-phase isothermal fluid flows. Other types of flows with nonlinear relations between  $v$, $\rho$ and $\nabla p$ can also be considered, see e.g. \cite{MuskatBook,BearBook,HK3,CHK4}.
\end{enumerate}
\end{remark}

\medskip
\noindent\textbf{Data availability.} 
No new data were created or analyzed in this study.

\medskip
\noindent\textbf{Funding.} A.I. obligatorily acknowledges OGRI's grant 122022800272-4.

\medskip
\noindent\textbf{Conflict of interest.}
There are no conflicts of interests.

\bibliography{paperbaseall}{}
\bibliographystyle{abbrv}
\end{document}